\documentclass[11pt]{amsart}
\usepackage{float}
\usepackage{amsfonts, amstext, amsmath, amsthm, amscd, amssymb, upgreek}
\usepackage{graphicx, color,  subfigure, wrapfig, overpic}
\usepackage[all,cmtip]{xy} 
\usepackage{tikz}

\let\oldmarginpar\marginpar
\renewcommand\marginpar[1]{\oldmarginpar[\raggedleft\footnotesize #1]%
{\raggedright\footnotesize #1}}

\AtBeginDocument{
   \def\MR#1{}
}
\newcommand{\Sp}{\rm{S}}
\newcommand{\C}{\mathbb{C}}
\newcommand{\R}{\mathbb{R}}

\newcommand{\Z}{\mathbb{Z}}

\newcommand{\HH}{\mathbb{H}}

\newcommand{\W}{\mathcal W}
\newcommand{\vol}{{\rm vol}}

\newcommand{\voct}{{v_{\rm oct}}}
\newcommand{\vtet}{{v_{\rm tet}}}
\renewcommand{\L}{\mathcal L}

\theoremstyle{plain}
\newtheorem{theorem}{Theorem}[section]
\newtheorem{corollary}[theorem]{Corollary}
\newtheorem{lemma}[theorem]{Lemma}
\newtheorem{prop}[theorem]{Proposition}

\newtheorem{conjecture}[theorem]{Conjecture}
\newtheorem{example}[theorem]{Example}

\newtheorem*{namedtheorem}{\theoremname}
\newcommand{\theoremname}{testing}

\theoremstyle{definition}
\newtheorem{define}[theorem]{Definition}

\newtheorem{remark}[theorem]{Remark}

\title [Fully Augmented Links in the Thickened Torus] {Fully Augmented Links in the Thickened Torus}

\author[Alice Kwon]{Alice Kwon}

\begin{document}
\maketitle

\begin{abstract}
  In this article we study the geometry of fully augmented link
  complements in the thickened torus and describe their
  geometric properties, generalizing the study of fully augmented
  links in $\Sp^3$. We classify which fully augmented links in the
  thickened torus are hyperbolic, show that their complements in the
  thickened torus decompose into ideal right-angled torihedra. 
  We also study volume
  density of fully augmented links in $\Sp^3$, defined to be the
  ratio of its volume and the number of augmentations.  We prove the Volume
  Density Conjecture for fully augmented links which states that the
  volume density of a sequence of fully augmented links in $\Sp^3$
  which diagrammatically converge, as defined below, 
  to a biperiodic link, converges to
  the volume density of that biperiodic link. Furthermore, we show that the complement of a
  sequence of these links approaches the complement of the biperiodic link as a geometric limit. 

\end{abstract}

\section{Introduction}
In this paper we study a class of links called \emph{fully augmented links}.
Fully augmented links in $\Sp^3$ are obtained from diagrams of
links in $\Sp^3$ as follows. Let $K$ be a link in $\Sp^3$
with a given planar link diagram $D(K)$. 
We encircle each 
twist region (a maximal string of bigons) of $D(K)$ with a single 
unknotted component called a \emph{crossing circle}.  The complement
of the resulting link is homeomorphic to the link obtained by removing 
all \emph{full-twists} i.e. pairs of crossings from each twist-region. 
Therefore, a diagram of the fully augmented link 
contains a finite number of crossing circles, each encircling two strands 
of the link. These crossing circles are perpendicular to the projection plane
and the other link components are embedded on the projection plane, except 
possibly for a finite number of single crossings, called 
\emph{half-twists} which are adjacent to the crossing circles. See Figure \ref{fig:falS3}.

The geometry of
fully augmented link complements in $\Sp^3$ can be explicitly described
in terms of an ideal right-angled polyhedral decomposition which is
closely related to the link diagram. This geometry has been studied in
detail by Adams \cite{CA}, Agol and D. Thurston \cite{Lackenby}, Purcell
\cite{JessicaP} and Cheesbro-Deblois-Wilton \cite{CDW}. In \cite{CKP2}
Champanerkar-Kofman-Purcell studied the geometry of alternating link
complements in the thickened torus and described their decompositions
into torihedra, which are toroidal analogs of polyhedra. In this paper we
combine the methods used to study fully augmented links in $\Sp^3$
and alternating links in the thickened torus to study the geometry of
fully augmented link complements in the thickened torus. We generalize
many geometric properties of fully augmented links in $\Sp^3$ to
those in the thickened torus, $T^2 \times I$, where $I = (-1,1)$.

 A \emph{biperiodic link} $\mathcal{L}$ is an infinite
 link in $\R^2 \times I$ with a projection on $\R^2 \times \{0\}$
 which is invariant under an action of a two dimensional lattice
 $\Lambda$ by translations. The quotient $L=\mathcal{L}/\Lambda$ is a
 link in $T^2 \times I$ with
 a projection on $T^2 \times \{0\}$. This projection on
 $T^2 \times \{0\}$ is the link diagram of $L$.

{\it Volume density} of a link $K$ was first introduced by Champanerkar,
Kofman and Purcell in \cite{CKPgmax} as the ratio of its hyperbolic volume, $\vol(K)$ and
its crossing number, $c(K)$. In \cite{CKPgmax} and \cite{CKP2} they studied
volume densities of sequences of alternating links in $\Sp^3$ which
diagrammatically converge to two specific biperiodic links called the
square weave and the triaxial link. They proved that volume density of
such a sequence of alternating links converges to that of the
corresponding biperiodic link. In general, they conjectured the
following:

\begin{conjecture} (Volume Density Conjecture
  \cite{CKP2}) \label{con:volumeDensityConjecture} Let $\mathcal{L}$
  be any biperiodic alternating link with alternating quotient link
  $L$. Let $\{K_n\}$ be a sequence of alternating hyperbolic links
  which F\o lner converges to $\mathcal{L}$. Then
  $$\displaystyle {\lim_{n \to \infty}} \frac{ \vol(K_n) }{c(K_n)} =
  \frac{\vol((T^2 \times I) - L)}{c(L)}.$$
\end{conjecture}

\begin{define}
A \emph{fully augmented biperiodic link} $\mathcal{L}$ is a fully augmented 
infinite link in $\R^2 \times I$ with a projection on $\R^2 \times \{0\}$
 which is invariant under an action of a two dimensional lattice
 $\Lambda$ by translations. The quotient $L=\mathcal{L}/\Lambda$ is a
fully augmented link in $T^2 \times I$ with
 a projection on $T^2 \times \{0\}$. 
\end{define}

We define the \emph{volume density} of a fully augmented link in $\Sp^3$ (with or without half-twists) to be the ratio of its volume and the number of augmentations. We similarly define volume density of fully augmented links in the thickened torus.
Using the geometry of fully augmented link complements in $\Sp^3$
studied previously and our results on the geometry of fully augmented
link complements in the thickened torus, we prove the Volume Density
Conjecture for fully augmented links.

In Section \ref{sec:one} we classify hyperbolic fully augmented links in the thickened torus. 

{
\renewcommand{\thetheorem}{\ref{thm:fal}}
\begin{theorem}
Let $K$ be a link in $T^2 \times I$ with a weakly prime, twist-reduced cellular link diagram $D$. Let $L$ be a link obtained by fully augmenting $D$. 
Then $T^2 \times I - L$ decomposes into two isometric totally geodesic right-angled torihedra, and hence $L$ is hyperbolic.
\end{theorem}
}

\begin{remark}
Augmented link diagrams are link diagrams obtained by adding crossing circles to some of the twist sites of a given link diagram and are different from  fully augmented links. It was proved in \cite{kwontham} that augmented links in the thickened torus are hyperbolic.  A generalization 
to thickened surfaces was also proved in \cite{ColinAdamsNew}. Theorem \ref{thm:fal} above gives a much stronger result for fully augmented links as it describes the right-angled geometry of the complement and uses very different proof techniques than \cite{kwontham} and \cite{ColinAdamsNew}. The decomposition of the link $L$ in Theorem \ref{thm:fal} into right-angled torihedra (see Definition \ref{def:angled-torihedron}) is very important for Theorem \ref{thm:volDetConjFAL} which investigates limit points of volume densities of fully augmented links. 
\end{remark}

 In Section \ref{sec:voldenconvconj} we discuss volume density and the volume density spectrum of fully augmented links in $\Sp^3$ and give many examples. In Section \ref{subsec:Folner} we define F\o lner convergence for fully augmented links and prove the volume density conjecture for fully augmented links. F\o lner convergence for links was first defined
in \cite{CKPgmax} for alternating links, we adapt the definition of F\o lner convergence in this paper for
sequences of fully augmented links.

{
\renewcommand{\thetheorem}{\ref{thm:volDetConjFAL}}
\begin{theorem} 
Let $\mathcal{L}$ be a biperiodic fully augmented link with quotient link $L$. Let $\{K_n\}$ be a sequence of hyperbolic fully augmented links in $\Sp^3$ such that $K_n$ F\o lner converges to $\mathcal{L}$ geometrically. Then $$\displaystyle {\lim_{n \to \infty}}  \frac{ \vol(K_n) }{a(K_n)} = \frac{\vol((T^2 \times I) - L)}{a(L)},$$ where $a(K)$ denoted the number of augmentations of a fully augmented link $K$. 
\end{theorem}  
}

As an application in Corollary \ref{cor:falAsymptoticVolume} we show that the end point $10\vtet$ of the volume density spectrum of fully
augmented links in $\Sp^3$ is a limit point, by constructing a sequence of hyperbolic fully augmented links in $\Sp^3$
which F\o lner converge everywhere to a fully augmented biperiodic link whose volume density is $10\vtet$.    

\ \\
{\bf Acknowledgements} I would like to thank my advisor Abhijit Champanerkar for guidance in this paper. I would also like to thank Ilya Kofman and Jessica Purcell for helpful conversations in regards to this project.

\section{Hyperbolicity of Fully Augmented Links in the Thickened Torus and Volume Bounds} \label{sec:one}

To define fully augmented links in the thickened torus we first need to define twist-reduced diagrams for links in $T^2 \times I$. Howie-Purcell defined twist-reduced diagrams for links in thickened surfaces in \cite{JessJosh}. However for links in the thickened torus we can also define twist-reduced diagrams using the biperiodic link diagram in $\R^2$:   
 
 \begin{figure}
 \centering
 \includegraphics [height=4cm]{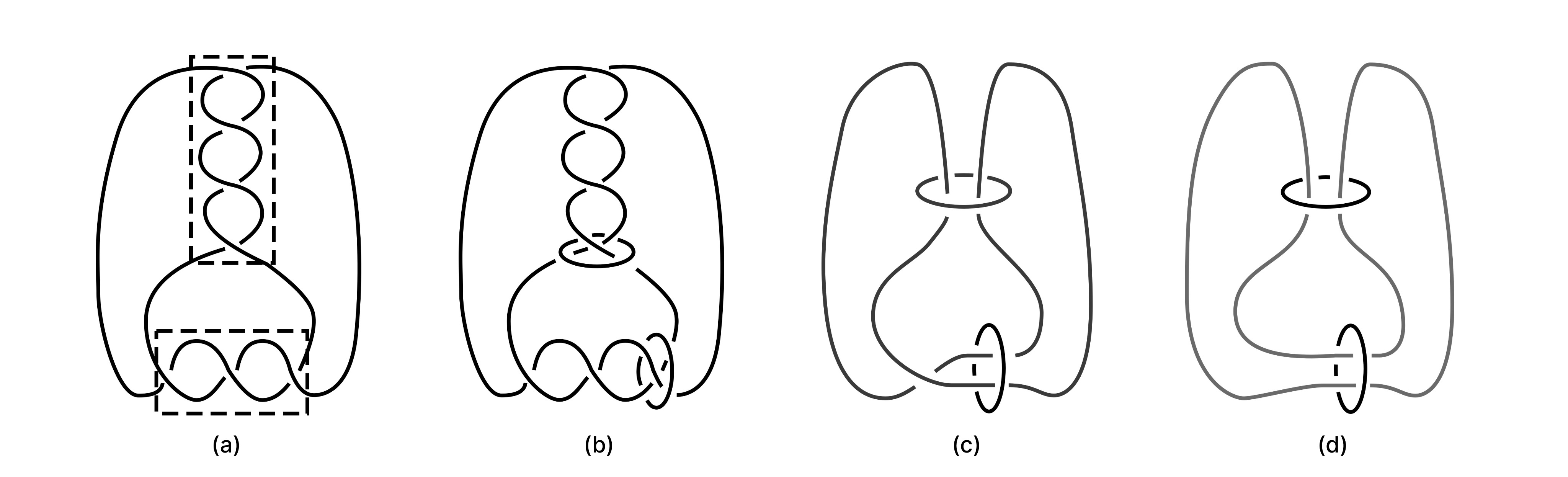}
 \caption{(a) Link diagram of $K$ (b) crossing circles added to each twist region (c) third picture is a fully augmented link diagram with all full-twists removed (d) fully augmented link diagram with no half-twists}
 \label{fig:falS3}
 \end{figure}


 \begin{define}
   A \emph{twist region} in the biperiodic link diagram $\mathcal{L}$
   is a maximal string of bigons, or a single crossing. A \emph{twist region}
   in the link diagram $L=\mathcal{L}/\Lambda$ is a quotient of a
   twist region in $\mathcal{L}$.
\end{define}

\begin{figure}
\includegraphics[height=3cm]{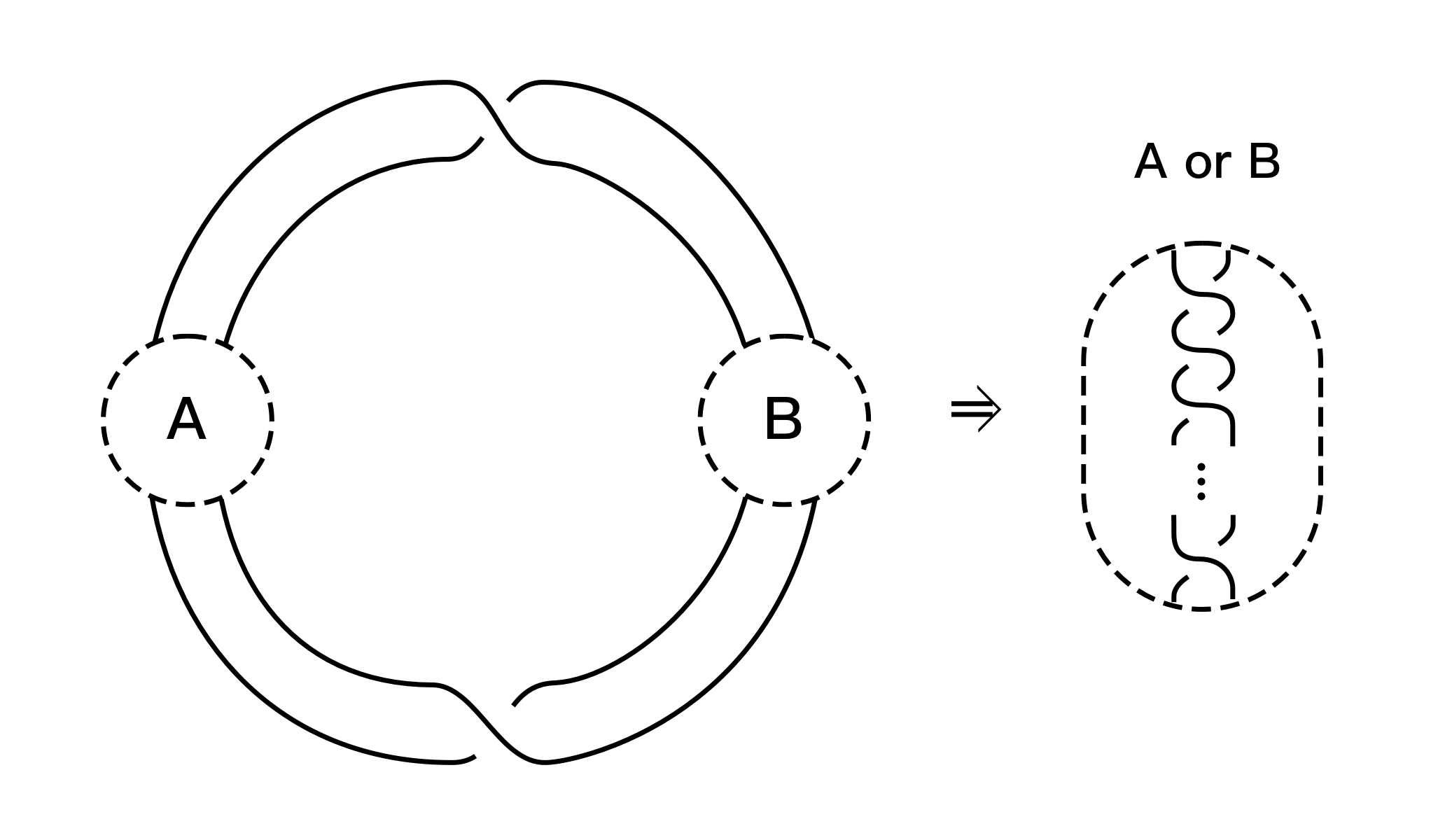}
\caption{Twist-reduced diagram}
\label{fig:twist-reduced}
\end{figure}

A biperiodic link $\mathcal{L}$ is called \emph{twist-reduced} if for
any simple closed curve on the plane that intersects $\mathcal{L}$
transversely in four points, with two points adjacent to one crossing
and the other two points adjacent to another crossing, the simple
closed curve bounds a subdiagram consisting of a (possibly empty)
collection of bigons strung end to end between these crossings.  See
Figure \ref{fig:twist-reduced}. We say $L$ is \emph{twist-reduced} if
it is the quotient of a twist-reduced biperiodic link.

\begin{figure}
\centering
 \begin{tabular}{cc}
\includegraphics[height =4cm]{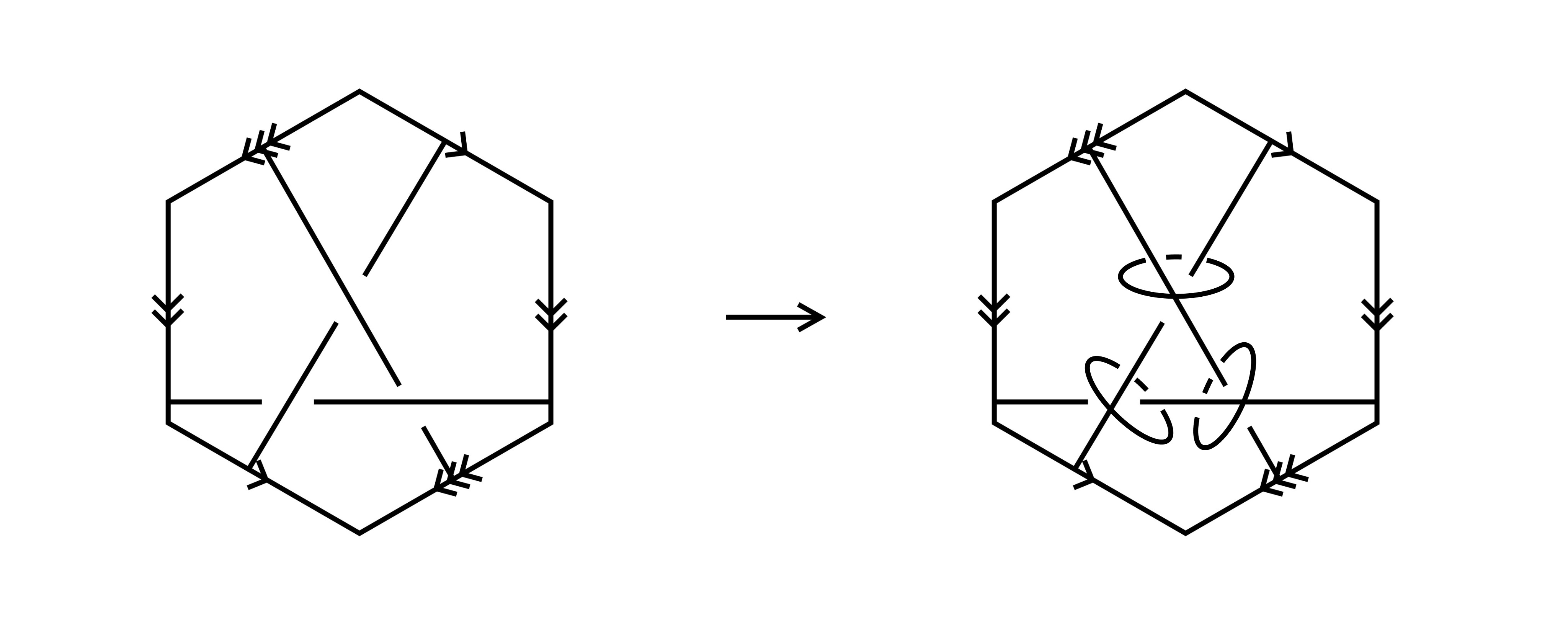}&\\
(a)\\
 \includegraphics [height=4cm]{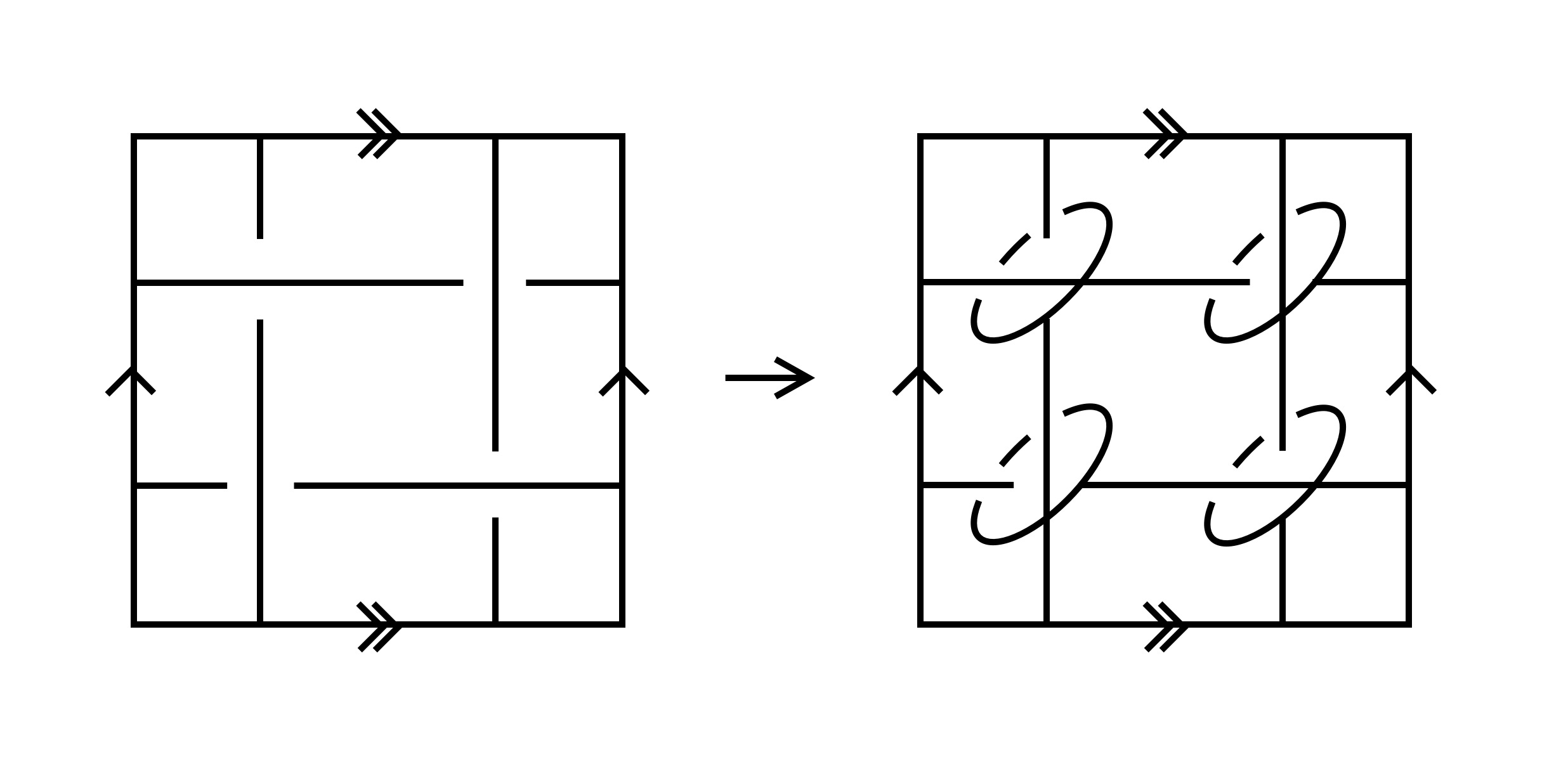}&\\
 (b)
  \end{tabular}
\caption{(a) Example of fully augmented triaxial link. (b) Example of fully augmented link on the square weave.}
\label{fig:fal}
\end{figure}

  \begin{define}\label{def:falinT2}
    A \emph{fully augmented link diagram in $T^2 \times I$} is a diagram of 
    a link $L$ that
    is obtained from a twist-reduced diagram $K$ in
    $T^2 \times I$ as follows: (1) augment every twist region with a
    circle component, called a \emph{crossing circle} (2) get rid of
    all full twists. See Figure \ref{fig:fal}. A {\it fully augmented
      link in $T^2 \times I$} is a link which has a fully augmented
    link diagram in $T^2 \times I$.
\end{define}

\begin{remark}
  For fully augmented links in $\Sp^3$, depending on the parity of the
  number of crossings in a twist region the fully augmented link may
  or may not have a half-twist at that crossing circle. See the third
  figure from the left in Figure \ref{fig:falS3}. Similarly depending
  on the parity of the number of crossings at a twist region a fully
  augmented link in the thickened torus may or may not have a
  half-twist at that crossing circle.
\end{remark}

\begin{define}
A graph $G = (V, E)$ on the torus is \emph{cellular} if its complement is a collection
of open disks.
\end{define}


Torihedra were first defined in \cite{CKP2} and play the role of polyhedra in polyhedral decompositions of link compliments in $\Sp^3$ e.g. it is proved in \cite{CKP2} that a complement of a link in the thickened torus decomposes into torihedra. Here we recall the definition of a torihedron:

\subsection{Torihedral decomposition}

\begin{define} \label{def:torihedron}
 A \emph{torihedron} is a cone on the torus, i.e. $T^2 \times [0,1]/(T^2 \times \{1\})$, with a cellular graph $G$ on $T^2 \times \{0\}$. The edges and faces of $G$ are called the edges and faces of the torihedron. An \emph{ideal torihedron} is a torihedron with the vertices of $G$ and the vertex $T^2 \times \{1\}$ removed.  Hence, an ideal torihedron is homeomorphic to $T^2 \times [0,1)$ with a finite set of points (ideal vertices) removed from $T^2 \times \{0\}$. The graph $G$ is called the \emph{graph of the torihedron}.
\end{define}

\begin{define} \label{def:angled-torihedron}
An \emph{angled torihedron} is a torihedron with an angle assignment on each edge of the graph of the torihedron. An assignment of $\pi/2$ angle on each edge is called a \emph{right-angled torihedron}. 
\end{define}

\begin{prop}\label{prop:Jess1}
Let $L$ be a fully augmented link in $T^2 \times I$, then there is a decomposition of the link complement $(T^2 \times I )- L$ into two combinatorially isomorphic torihedra such that: 
\begin{enumerate}
\item The faces of each torihedron can be checkerboard colored such
  that the shaded faces are triangular and arise from the bow-ties
  corresponding to crossing circles;
\item The graph of each torihedron is 4-valent.
\end{enumerate}
\end{prop}
 
{\it Proof.} We follow the cut-slice-flatten construction described in
\cite{Lackenby}.  Let $L$ be a fully augmented link in
$T^2 \times I$. We begin by assuming that there are no half twists,
the crossing circles are lateral to $T^2 \times \{0\}$, and the
components of $L$ that are not crossing circles lie flat on
$T^2 \times \{0\}$. There are twice punctured disks
bounded by the crossing circles which are perpendicular to the
projection plane.

\begin{figure}
\centering
\begin{tabular}{cccc}
\includegraphics [height=3cm]{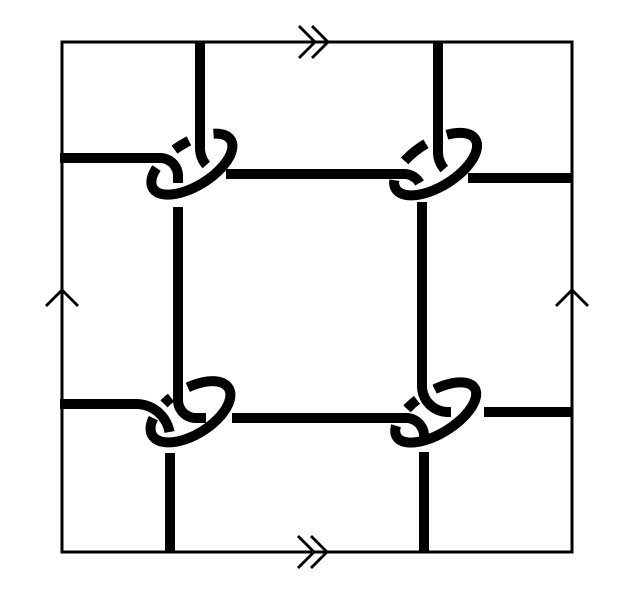}&
\includegraphics [height=3cm]{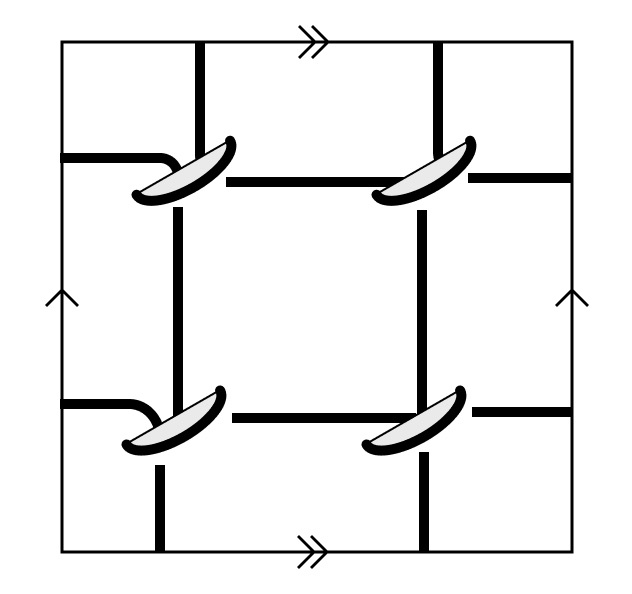}&
\includegraphics [height=3cm]{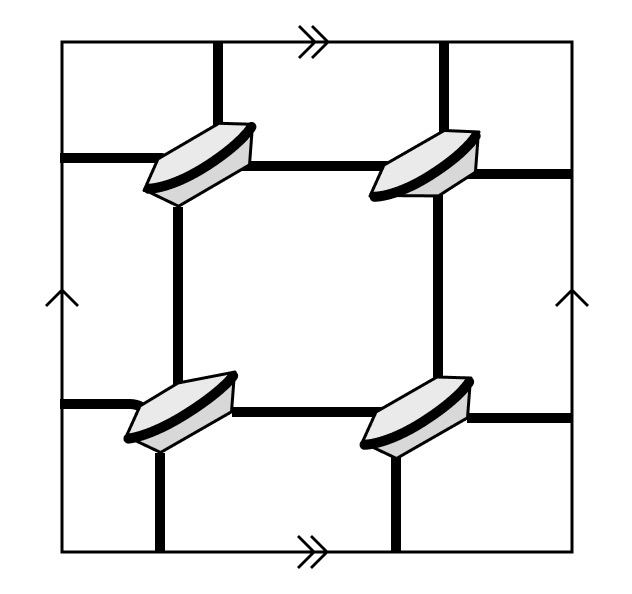}&
\includegraphics [height=3cm]{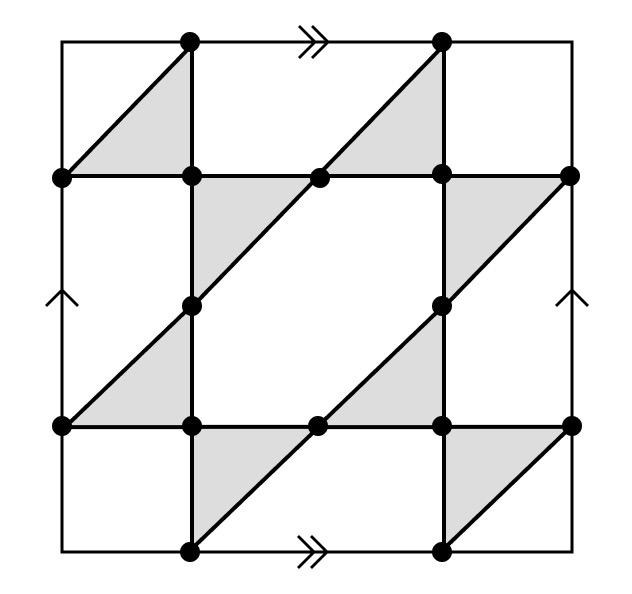}\\
(a)&(b)&(c)&(d)
\end{tabular}
\caption{(a) A fundamental domain for a fully augmented square weave, $L$. (b) Disks cut in half at each crossing circle. (c) Sliced and flattened half-disks at each crossing circle  (d) Collapsing the strands of the link and parts of the augmented circles (shown in bold) to ideal points gives the bow-tie graph $\Gamma_{L}$. The disks become shaded bow-ties and the white regions become hexagons.}

\label{fig:falDecomp}
\end{figure}

\begin{enumerate}
\item Cut $T^2 \times I$ along the projection surface 
$T^2 \times \{0\}$ into two pieces. 
This cuts each of the twice-punctured disks bounded by a 
crossing circle in half. See Figure \ref{fig:falDecomp}(b). 
\item For each of the two pieces resulting from Step 1, slice the
  middle of the halves of twice-punctured disks and flatten the half
  disks out. See Figure \ref{fig:falDecomp}(c).
\item Collapse strands of the link and parts of the augmented circles
 to ideal vertices in each of the two pieces. See Figure \ref{fig:falDecomp}(d). 
\end{enumerate}

It follows from Steps (1), (2) and (3) that each piece
of the decomposition is homeomorphic to $T^2 \times [0,1)$ 
with the same graph on $T^2 \times \{0\}$ with vertices deleted.
Hence $(T^2 \times I) - L$
decomposes into two identical ideal torihedra. 

After Step (2) the cut-sliced-flattened half-disks become a hexagon 
with an edge in the middle corresponding to the strand of the
half of a crossing circle. Upon collapsing the crossing circle 
this becomes a bow-tie. See Figures \ref{fig:falDecomp}(c) and
\ref{fig:falDecomp}(d). Each vertex of the graph is four-valent since
it is shared by two triangles of either two different bow-ties or one
bow-tie.  Again by construction each edge is shared by a triangle of a
bow-tie and a polygon that does not come from a bow-tie see Figure
\ref{fig:falDecomp}(d). Hence we can shade each triangle of the
bow-tie to get a checkerboard coloring on the graph of the torihedron
such that the shaded faces are bow-ties. 

\begin{figure}
 \centering
 \begin{tabular}{cc}
 \includegraphics [width=9cm]{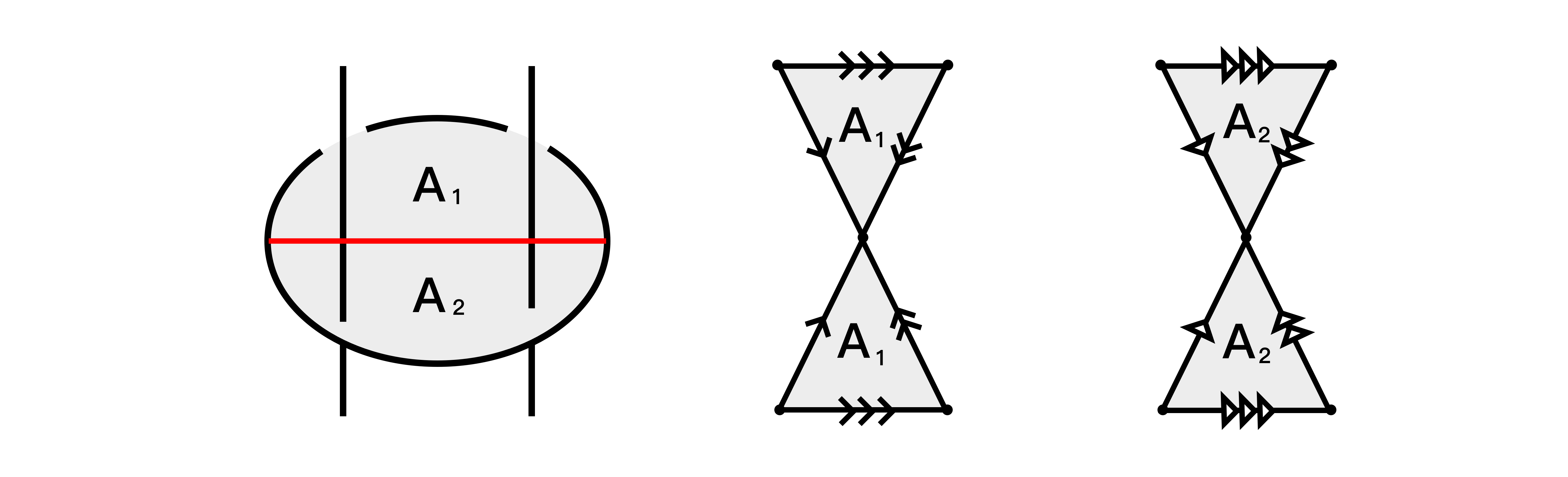}&
 \includegraphics [width=7cm]{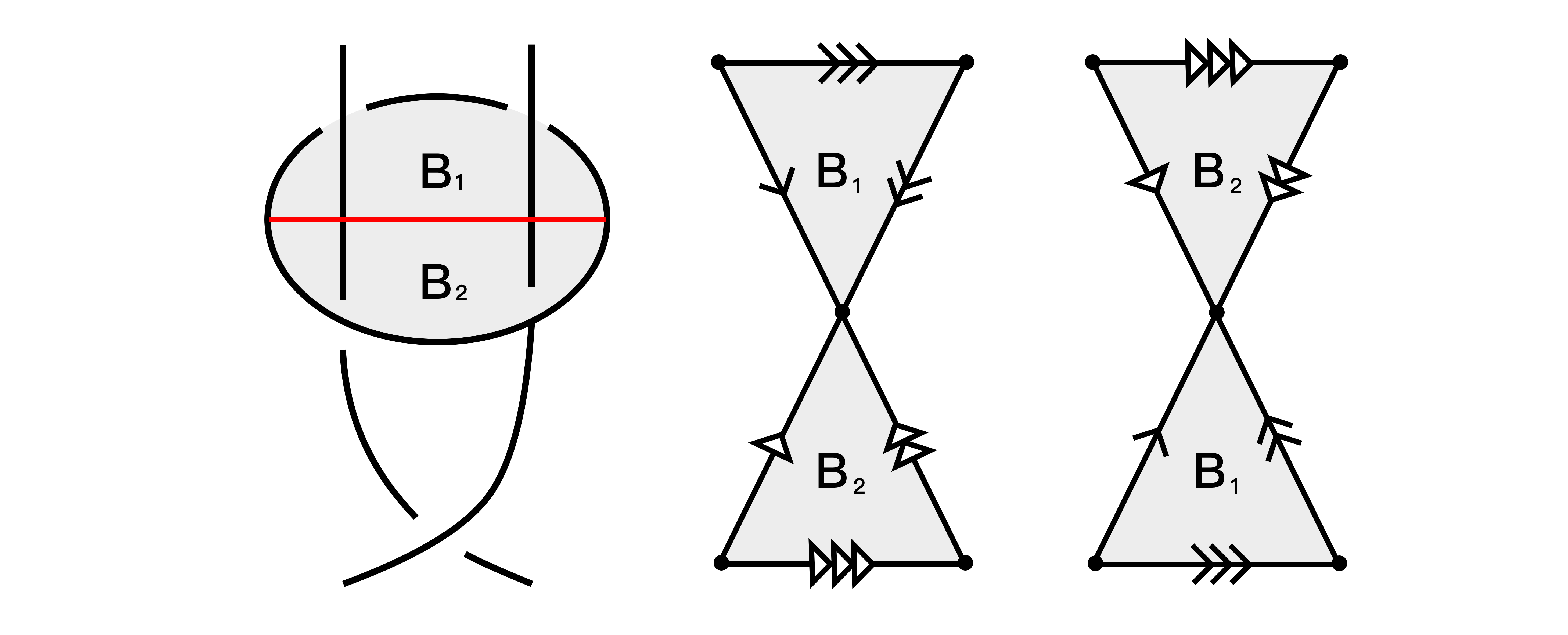}\\
 (a)&(b)
 \end{tabular}
 \caption{(a) Gluing information on the edges of the bow-tie without half-twists. (b) Gluing information on the edges of the bow-tie with half-twist.}
 \label{fig:falGluings}
 \end{figure}

The two torihedra are glued together as follows: the white faces are glued 
to the corresponding white faces, and the bow-ties are glued as shown 
in Figure \ref{fig:falGluings}(a).  

\begin{figure}
 \centering
 \includegraphics[width=16cm]{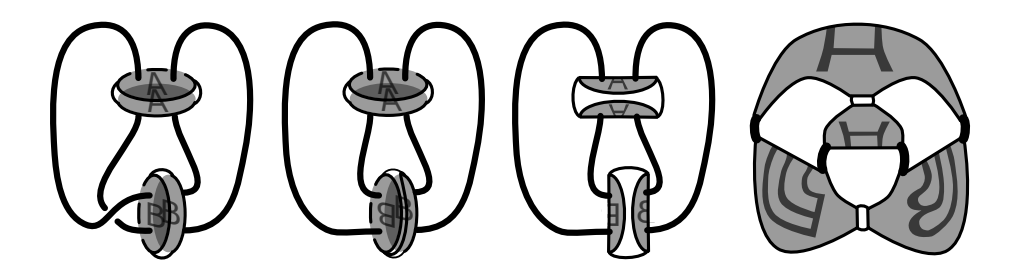}
 \caption{The gluing of the torihedra
 when half-twist is present (disk B) and when half-twist is absent (disk A). 
 This figure taken from \cite{JP-Book} is for links in $S^3$, but since this is a local move, the same gluing works for links in $T^2 \times I$.}
 \label{fig:falGluings1}
 \end{figure}

In the case when there is a half-twist at a crossing circle, we split the whole twice-punctured disk into two copies, and flip one of the disks to remove the half-twist. This only affects the gluing of faces of the torihedra. 
Hence if there are half-twists, then we get the same torihedra but with a
different gluing pattern on the bow-ties as shown in Figures
\ref{fig:falGluings1} and \ref{fig:falGluings}(b). \qed

\begin{define}\label{def:bowtieGraph}
For a fully augmented link $L$ in the thickened torus the decomposition of $T^2 \times I - L$ described above is called the \emph{bow-tie torihedral decomposition} of $L$. We call the graph of the torihedra the \emph{bow-tie graph} of $L$ and denote it as $\Gamma_{L}$. 
\end{define}

\begin{lemma} \label{lem:Jess}
Let $L$ be a hyperbolic fully augmented link in $T^2 \times I$. The following surfaces are embedded totally geodesic surfaces in the hyperbolic structure on the link complement. 
\begin{enumerate}
\item Each twice-punctured disk bounded by a crossing circle,
\item Each connected component of the projection surface. 
\end{enumerate}
\end{lemma}

{\it Proof.} (1) The disk $E$ bounded by a crossing circle is
punctured by two arcs of the link diagram lying on the projection
plane.
Adams \cite{adamsTP} showed that any incompressible
twice-punctured disk properly embedded in a hyperbolic 3-manifold is
totally geodesic.
 Hence it suffices to show that $E$ is incompressible.  Let $L$ be a
hyperbolic fully augmented link in $T^2 \times I$.  Since
$T^2 \times I \simeq \Sp^3-H$, where $H$ is the Hopf link, $L \cup H$
is a hyperbolic link in $\Sp^3$. 

Suppose there is a compressing disk $D$ with with
$\partial D \subset E$. Since $\partial D$ is an essential closed
curve on $E$, it must encircle one or two punctures of $E$.  Suppose
it encircles only one puncture. This means that the union of $D$ and
the disk bounded by $\partial D$ inside the closure of $E$ forms a
sphere in $\Sp^3$ met by the link exactly once. This is a
contradiction to the generalized Jordan curve theorem. Hence
$\partial D$ must bound a twice-punctured disk $E'$ on $E$. This means
$\overline{(E -E')} \cup D$ is a boundary compressing disk for the
crossing circle, contradicting the boundary irreducibility of 
${\rm S}^3 - (L \cup H)$.

\ \\
\indent (2) Notice that the reflection through the projection surface
($T^2 \times \{0\}$) preserves the link complement, fixing the plane
pointwise. Then it is a consequence of Mostow-Prasad rigidity that
such a surface must be totally geodesic. 
(See Lemma 2.1 in \cite{JessicaP}) \qed

\subsection{Hyperbolicity}
\begin{figure}
\includegraphics[height=4cm]{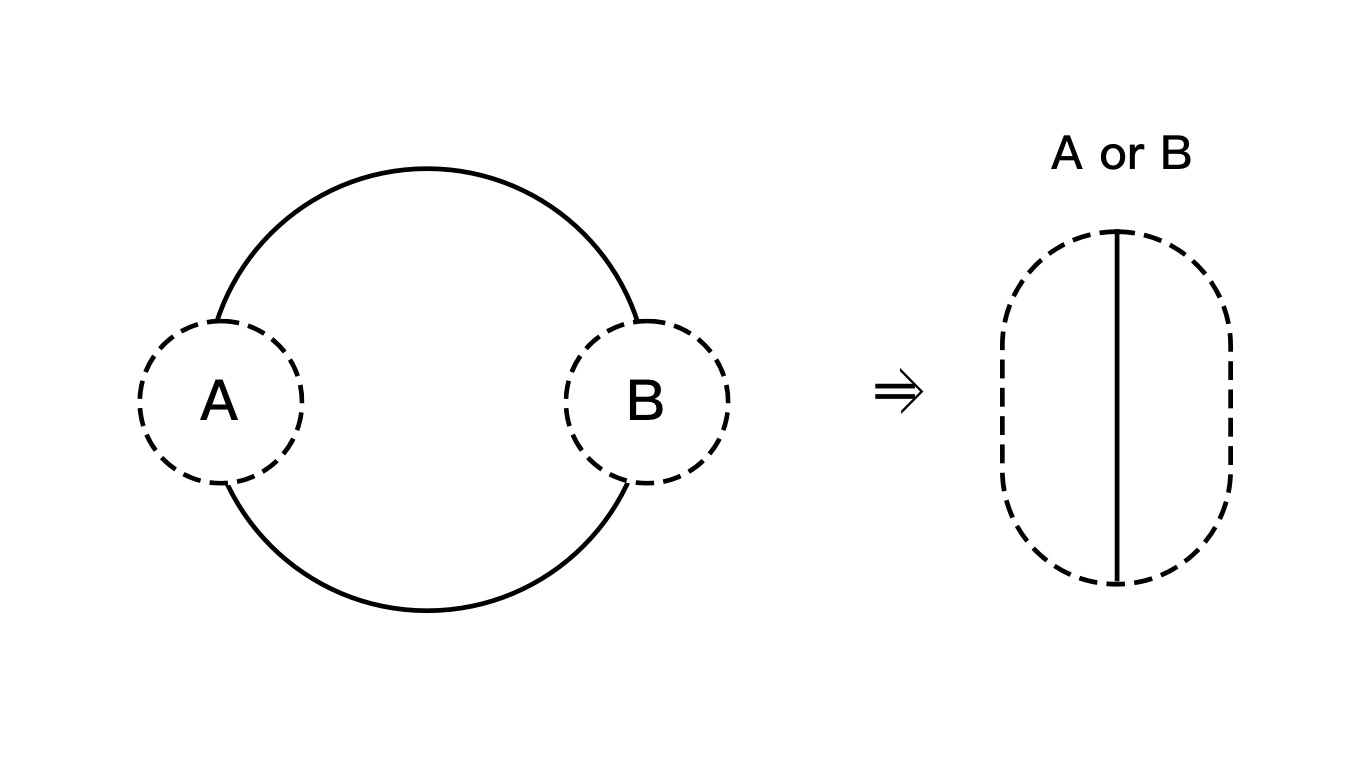}
\caption{Prime diagram}
\label{fig:prime}
\end{figure}

\begin{define} \label{def:weaklyprime}  Let $\L$ be a 
biperiodic link with diagram $D(\L)$. We say $D(\L)$ is prime if whenever  
a disk embedded in $\R^2 \times \{0\}$ meets $D(\L)$ transversely 
in exactly two edges, then the disk contains a simple edge of the diagram 
and no crossings. See Figure \ref{fig:prime}.

A diagram of a link $L$ in $T^2 \times I$, denoted $D(L)$ is \emph{weakly prime}
if $D(L)$ is a quotient of a prime biperiodic link diagram $D(\L)$ in $\R^2 \times \{0\}$ 
\end{define}

\begin{theorem}\label{thm:fal}
Let $K$ be a link in $T^2 \times I$ with a weakly prime, twist-reduced cellular link diagram $D$. Let $L$ be a link obtained by fully augmenting $D$. 
Then $T^2 \times I - L$ decomposes into two isometric totally geodesic right-angled torihedra, and hence $L$ is hyperbolic.
\end{theorem}

The proof of Theorem \ref{thm:fal} relies on a result about the existence of certain circle patterns on the torus due to Bobenko and Springborn \cite{BandS}. We use similar ideas from \cite{CKP2} to prove Theorem \ref{thm:fal} 

\begin{theorem} \label{thm:BandS}
\cite{BandS}. Suppose $G$ is a 4-valent graph on the torus $T^2$, and $\theta \in (0,2\pi)^E$ is a function on edges of $G$ that sums to $2\pi$ around each vertex. Let $G^*$ denote the dual graph of $G$. Then there exist a circle pattern on $T^2$ with circles circumscribing faces of $G$ (after isotopy of $G$) and having exterior intersection angles $\theta$, if and only if the following condition is satisfied:
\ \\
\indent Suppose we cut the torus along a subset of edges of $G^*$, obtaining one or more pieces. For any piece that is a disk, the sum of $\theta$ over the edges in its boundary must be at least $2\pi$, with equality if and only if the piece consists of only one face of $G^*$ (only one vertex of $G$).
\\
The circle pattern on the torus is uniquely determined up to similarity.
\end{theorem}

{\it Proof of Theorem \ref{thm:fal}}.  Decompose $(T^2 \times I) -L $
into two torihedra using Proposition \ref{prop:Jess1}. Let $\Gamma_L$
be the bow-tie graph on $T^2 \times \{0\}$. Assign angles
$\theta(e) = \pi/2$ for every edge $e$ in $\Gamma_L$. We now verify
the condition of Theorem \ref{thm:BandS}. This will prove the
existence of an orthogonal circle pattern 
(circle pattern whose angle at the intersection of any two circles is orthogonal) 
circumscribing the faces of
$\Gamma_L$.
\ \\
\indent Let $C$ be a loop of edges of $\Gamma_L^*$ enclosing a disk
$D$. Suppose $C$ intersects $n$ edges of $\Gamma_L$ transversely. Let $V$
denote the number of vertices of $\Gamma_L$ that lie in $D$, and let
$E$ denote the number of edges of $\Gamma_L$ inside $D$ disjoint from
$C$. Because the vertices of $\Gamma_L$ are 4-valent and since the
edges inside $D$ which are disjoint from $C$ get counted twice for
each of its end vertex, $n+ 2 E = 4V $. This implies $n$ is
even. Since $K$ is weakly prime and $C$ is made up of edges dual to
$\Gamma_L$ this implies $n > 2$. Since $n$ is even, $n \geq 4$. Hence
the sum of the angles for all edges of $C$ must be at least $2\pi$.
\ \\
\indent We now show that this is an equality if and only if $C$
consists of one face of $\Gamma_L^*$, i.e. $C$ encloses only one
vertex. Suppose that the sum
$\displaystyle{\sum_{e \in C}}\theta(e) > 2\pi$. Since
$\theta(e)= \pi/2$, for every $e \in \Gamma_L$, and $n$ is even,
$n \geq 6$.  Moreover
$$n \geq 6 \implies 4V-2E \geq 6 \implies 2V-E \geq 3 \implies V\geq
2.$$ Hence $C$ encloses more than one vertex.
\ \\
\indent Conversely, let
$\displaystyle{\sum_{e \in C}}\theta(e) = 2\pi$. This implies $n=4$.

Let the edges of $C$ be $e_i$ for $0 \leq i \leq 3$ with $e_i$ incident to vertices $v_{i}$ and $v_{i+1}$ and $v_0=v_4$. Let the faces dual to $v_i$ be $F_{v_i}$.  Without loss of generality, let $F_{v_0}$ be a shaded triangular face. Since $\Gamma_L$ is checkerboard colored, $F_{v_2}$ is also shaded triangular face.

\begin{figure}
\centering
\begin{tabular}{cc}
\includegraphics [height=7cm]{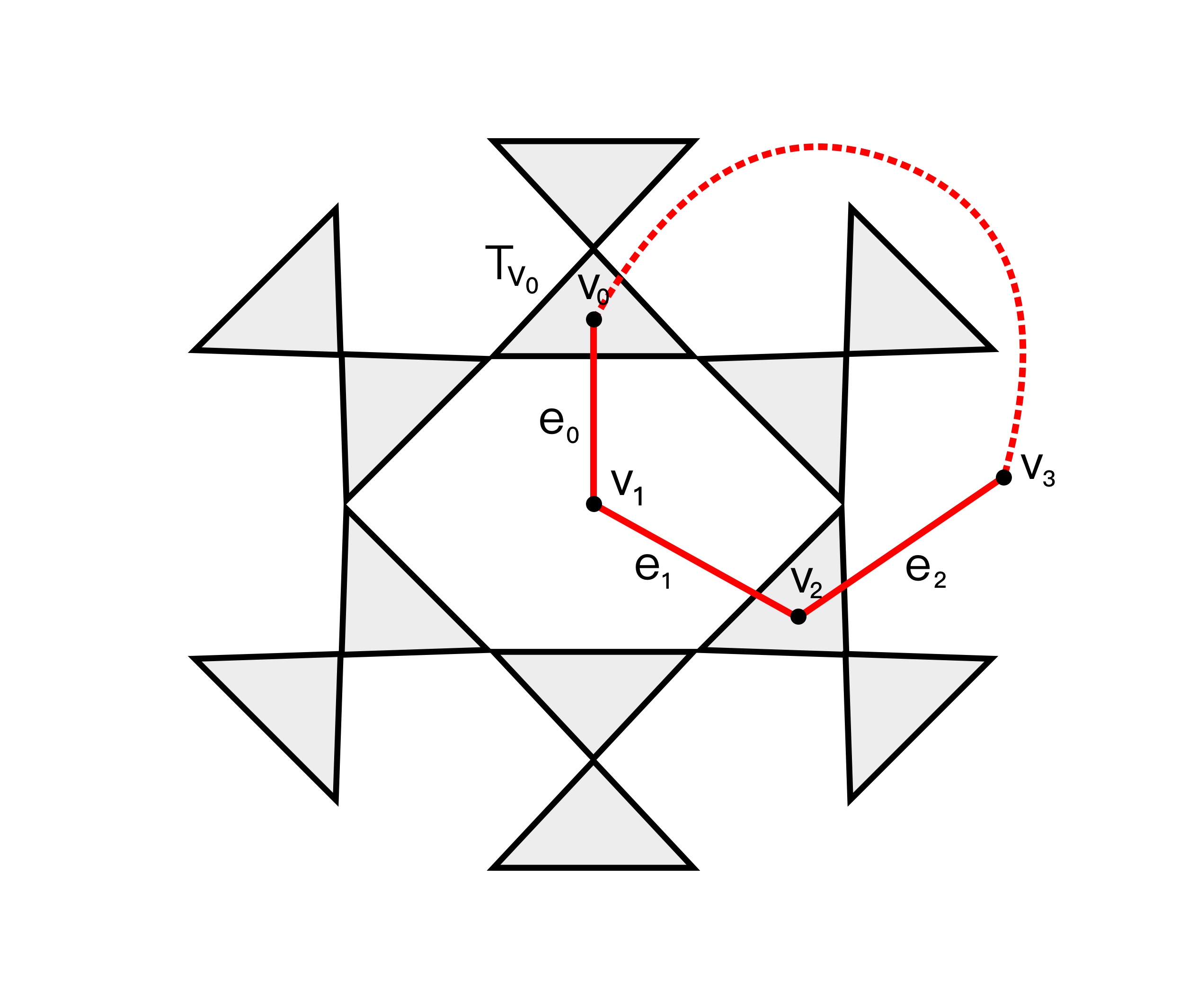}&
\includegraphics [height=7cm]{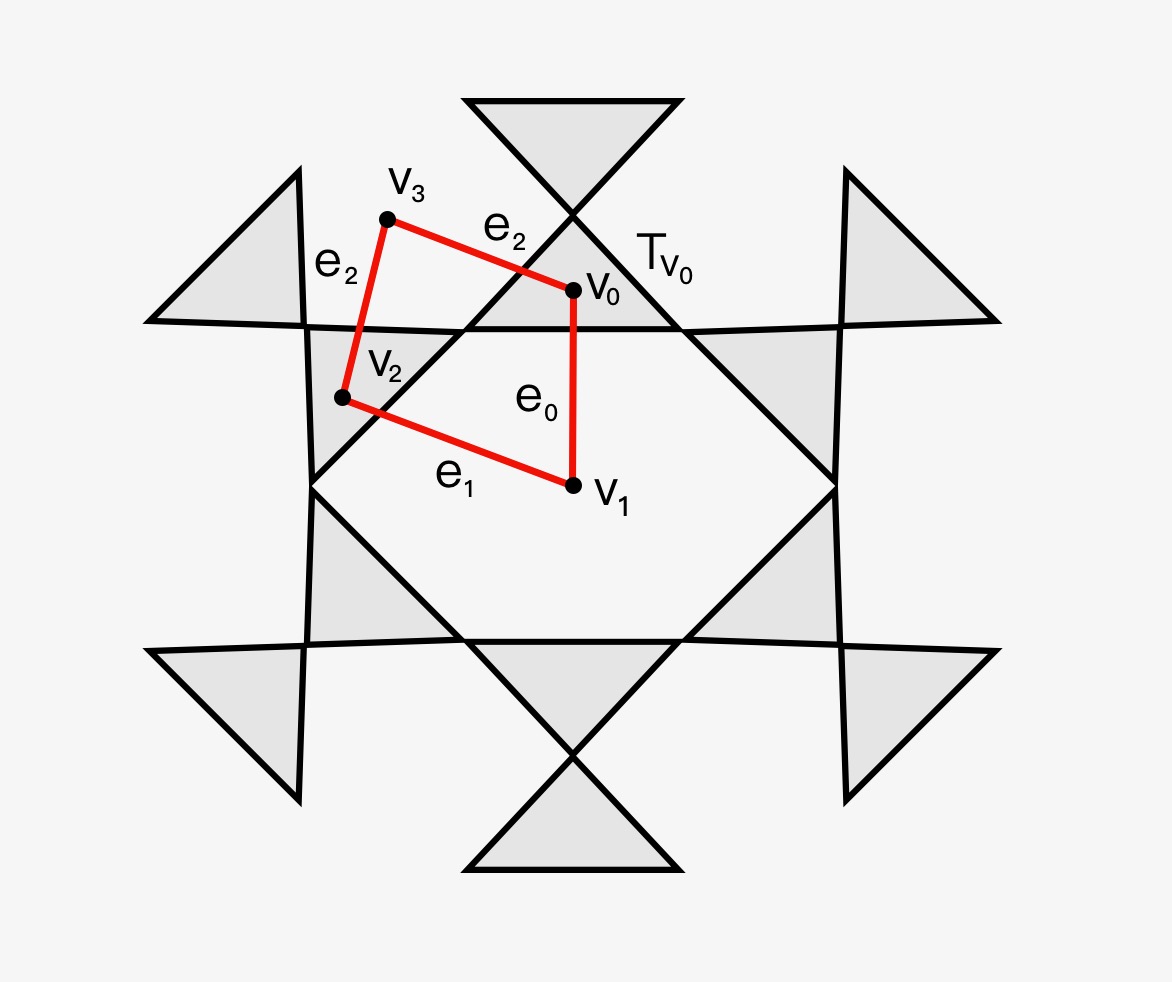}\\
\\
(a)&(b)
\end{tabular}
\caption{(a) When $n \geq 5$ and $C$ closes with $\geq 5$ edges. (b) When $n=4$ and $C$ closes with 4 edges.}
\label{fig:gammaL}
\end{figure}

Suppose $F_{v_0} \cap F_{v_2} = \emptyset$ then the edge $e_2$ must enter a white face $F_{v_3}$ which has empty intersection with $F_{v_0}$. See Figure \ref{fig:gammaL}(a). 

\begin{figure}
\centering
\begin{tabular}{ccc}
\includegraphics [height=3cm]{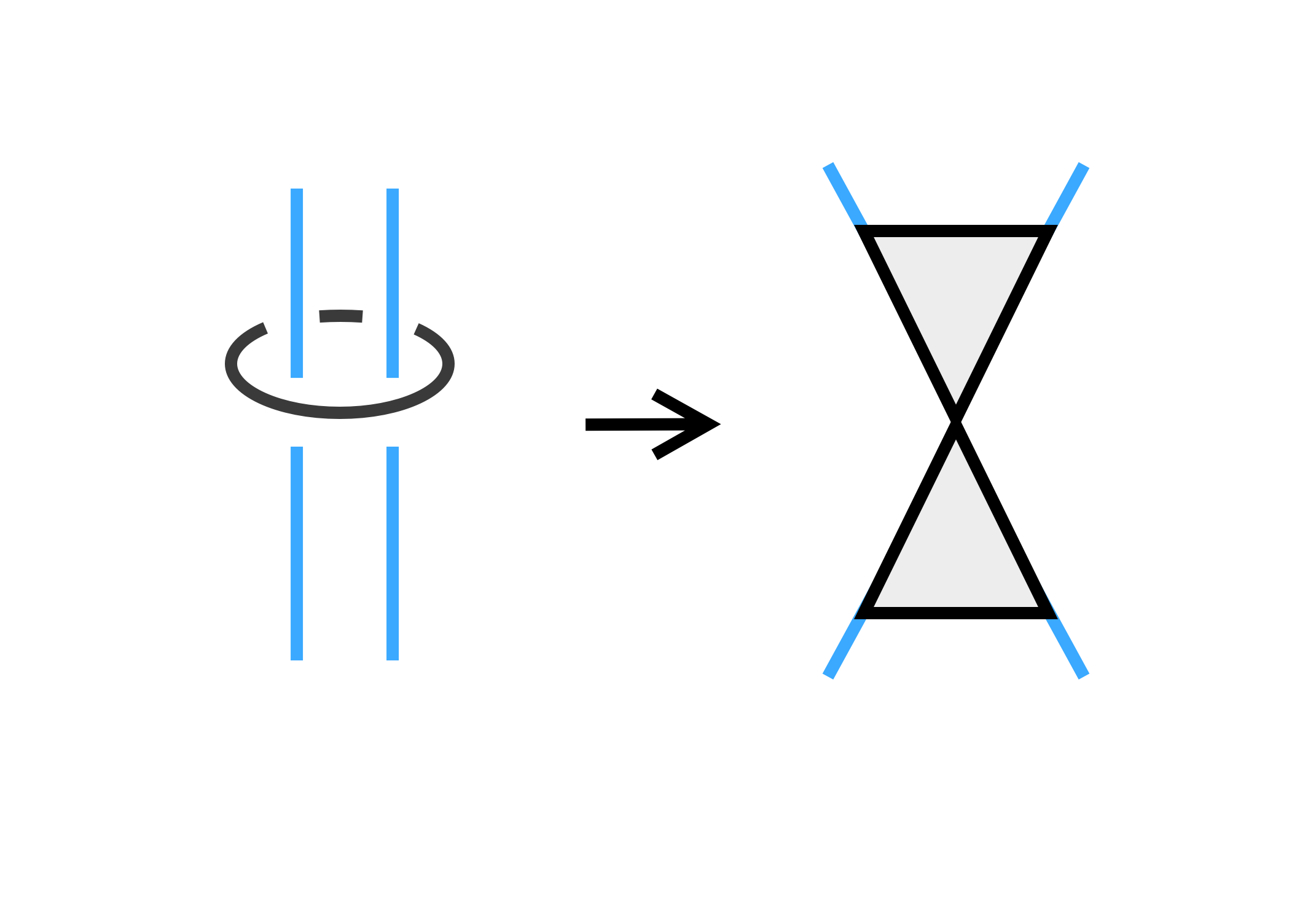}&
\includegraphics[height=9cm]{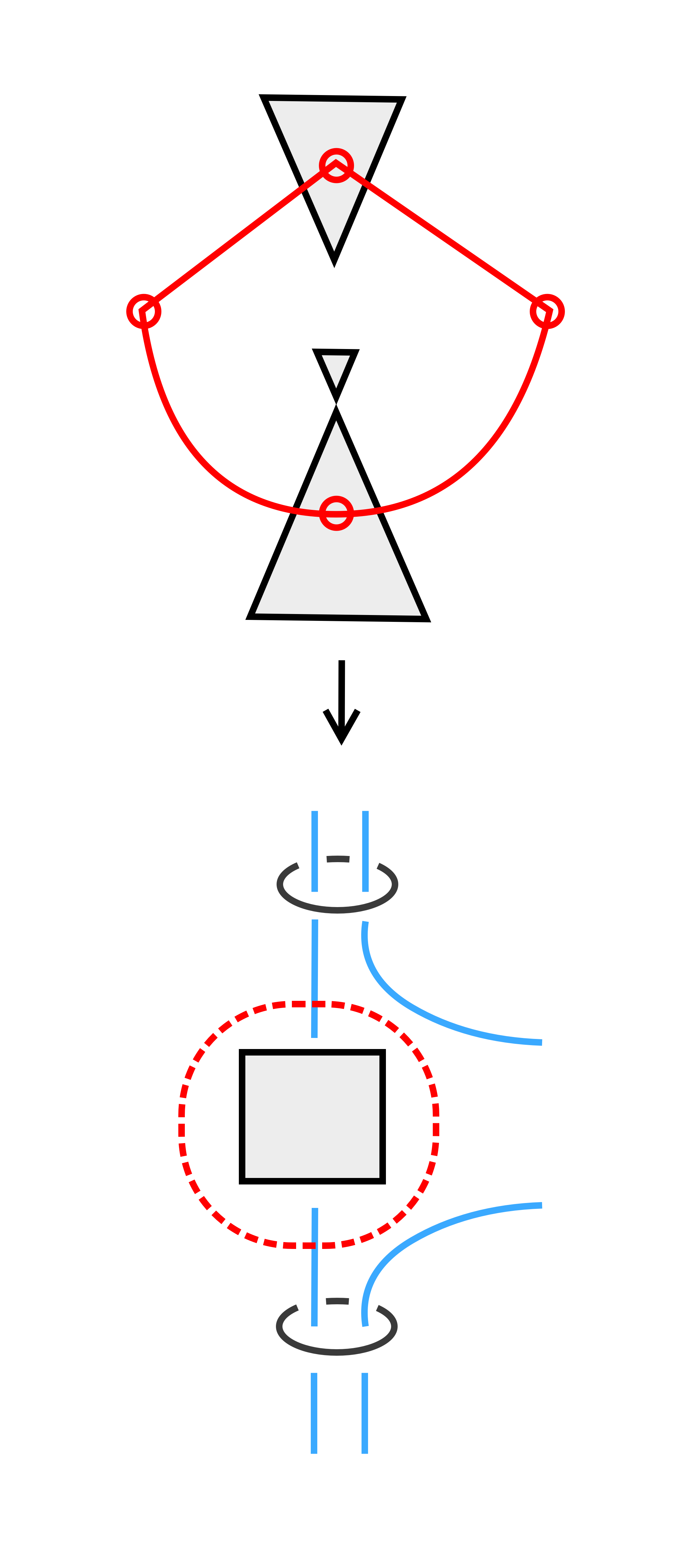}&
\includegraphics[height=9cm]{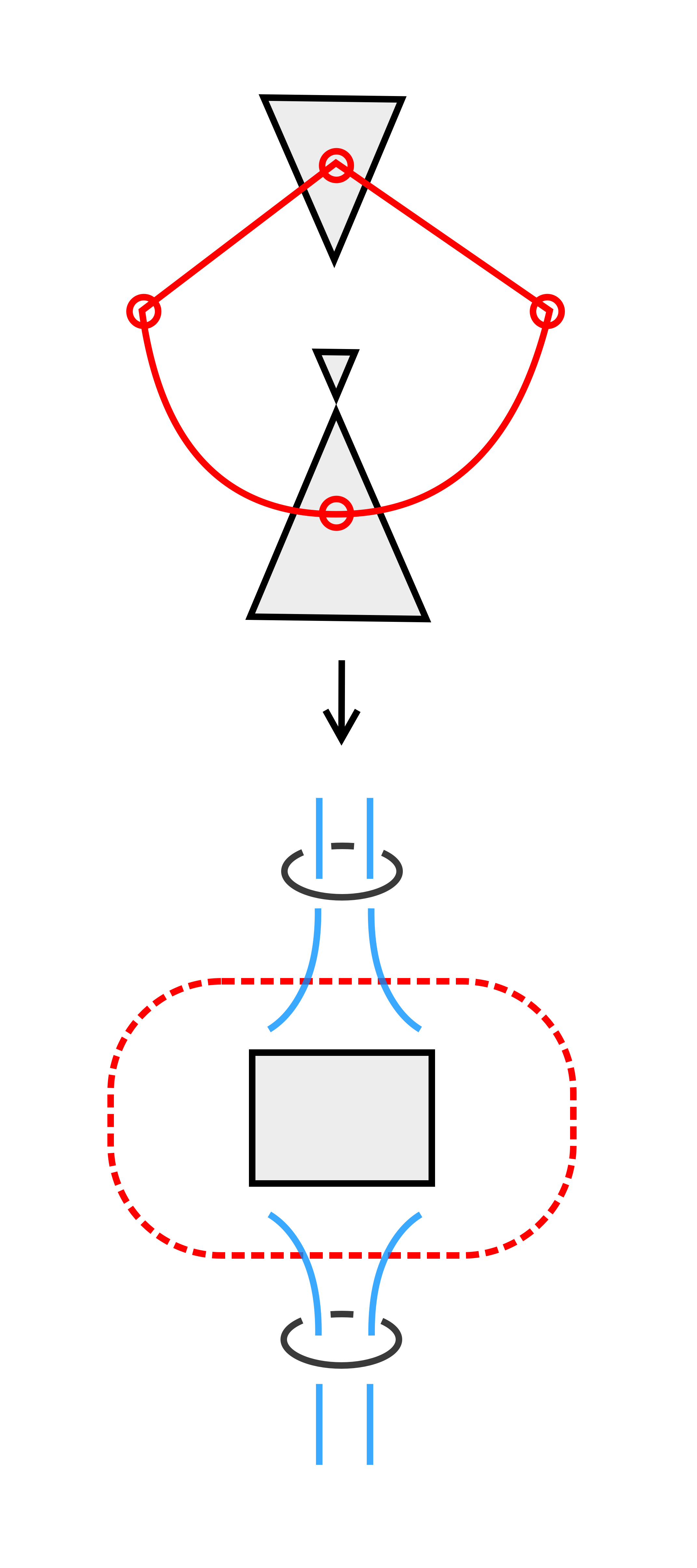}\\
(a)&(b)&(c)
\end{tabular}
\caption{(a) The crossing circle splits into a bowtie. (b) $C$ is in red, $C$ intersects the original link in two points and hence must be a trivial edge. (c) $C$ is in red, $C$ can intersect on the original link at four points and therefore must bound a twist region on one side.  }
\label{fig:gammaL0}
\end{figure}

Since the bowties correspond to crossing
circles, see Figure \ref{fig:gammaL0}(a), the loop $C$ gives a loop which
 intersects $L$. At the vertex $v_0$, which is in the shaded bowtie, at least one of the edges incident 
 to $v_0$ has to intersect $L$. If only one edge at $v_0$ intersects $L$, since $C$ bounds a disk, only one edge at $v_2$ intersects $L$, giving the case shown in Figure \ref{fig:gammaL0}(b). Similarly, If both edges incident to $v_0$ intersect $L$, since $C$ bounds a disk, then the same is true for both edges incident at $v_2$, giving the case shown in Figure \ref{fig:gammaL0}(c).
If $C$ intersects two strands of $L$ as in Figure \ref{fig:gammaL0}(b), since
$C$ bounds a disk, this contradicts the weakly prime condition of $K$.
If $C$ intersects two strands on each side as in Figure
\ref{fig:gammaL0}(c), this will contradict twist-reduced condition on
$K$. 

Therefore $F_{v_0} \cap F_{v_2} \neq \phi$. Since both faces are
triangles, they can only intersect in a vertex. This implies that $C$
encloses a single vertex. See Figure \ref{fig:gammaL}(b).



Now, since we showed that $\Gamma_L$ is a graph on the torus which
satisfies the conditions of Theorem \ref{thm:BandS} there exists an
orthogonal circle pattern on the torus with circles circumscribing the
faces of $\Gamma_L$. Since a white face of the decomposition intersects 
any other white face only at ideal
vertices, the circles which circumscribe the white faces create a
circle packing, where the points of tangency are those corresponding to the
associated ideal vertices. Since $\Gamma_L$ is 4-valent and every edge
has been assigned an angle of $\pi/2$, the circles of the
shaded faces meet orthogonally.

Lifting the circle pattern to the universal cover of the torus defines
an orthogonal biperiodic circle pattern on the plane. Considering the plane $z=0$ as a part 
of the boundary of $\mathbb{H}^3$, this circle pattern
defines a right-angled biperiodic ideal hyperbolic polyhedron in
$\mathbb{H}^3$. The torihedron of the decomposition of
$(T^2 \times I) - L$ is the quotient of $\HH^3$ by $\Z \times \Z$ which is now
realized as a right-angled hyperbolic torihedron. It follows from
Theorem 1.1 of \cite{JessJosh} that $(T^2 \times I) -L$ is hyperbolic.\qed


\begin{remark}
  Adams \cite{CA} proved that fully augmented link complements in $\Sp^3$ are
  hyperbolic. We have proved an analogous result for fully augmented
  link complements in $T^2\times I$.  Our method of finding an
  orthogonal circle pattern which circumscribed the faces of the bow-tie
  graph can also be applied to the case of fully augmented links in $\Sp^3$. In
  this we have to use Andreev's theorem \cite{Thurston} to ensure a totally
  geodesic right-angled polyhedra.

\end{remark}

 \subsection{Volume bounds}
 We show that a hyperbolic fully augmented link with $c$ crossings in the thickened torus have an upper volume bound
 of $10c\vtet$. In the next section we show volume density convergence of fully augmented links. This means if we 
 can find a link in the thickened torus whose volume is exactly $10c\vtet$ the corresponding biperiodic link
 will have volume density $10\vtet$. We will use this to show that an end point of the volume density spectrum 
 of fully augmented links can be obtained as a limit.  
 \begin{prop} \label{prop:lowerBoundFal}
 Let $L$ be a hyperbolic fully augmented link with $c$ crossing circles. Then $$2c\voct  \leq {\rm vol}(T^2 \times I - L) \leq 10c\vtet$$
where $\voct = 3.66386...$ is the volume of a regular ideal octahedron and $\vtet = 1.01494...$ is the volume of a regular ideal tetrahedron. 
 \end{prop}

{\it Proof.} We will first prove the lower bound. By work of Adams \cite{CA},
 the volume of the complement of $L$ in $T^2 \times I$
 agrees with that of the fully augmented
 link with no half-twists. This means a lower volume bound for the complement of $L$
 in $T^2 \times I$ 
 with half-twists will be a lower volume bound of the fully augmented link with no half-twists. 
 Hence we will assume $L$ has no half-twists and obtain a lower bound for $T^2 \times I - L$. 
 
 Cut $T^2 \times I - L$ along the reflection plane $T^2 \times \{0\}$,
 dividing it into two isometric hyperbolic manifold. The boundary of
 each of these are the regions of $L$ on the projection surface with
 punctures for the crossing circles. By Lemma \ref{lem:Jess} these regions 
 are geodesic. Hence cutting along the projection surface divides
 $T^2 \times I - L$ into isometric hyperbolic manifolds with totally
 geodesic boundary.

\indent Miyamoto showed that if $N$ is a hyperbolic 3-manifold with totally
 geodesic boundary, then vol($N$) $\geq -\voct\chi(N)$ \cite{Miya},
 with equality exactly when $N$ decomposes into regular ideal
 octahedra. In our case, the manifold $N$ consists of two copies of
 $T^2 \times [0,1)$ with half-annuli removed for half the crossing circles.
 For every half a crossing circle removed, we are removing one edge
 and two vertices. Hence for each crossing circle removed the Euler
 characteristic changes by $-1$. Since there are $c$ crossing circles the
 Euler characteristic would be $-c$ for each half-cut
 $T^2 \times [0,1)$. The lower bound now follows. 
 
 We now prove the upper bound.  The torihedral decomposition of the
 link complement gives a decomposition into two identical ideal
 torihedra. Every triangular shaded face which comes from a bowtie
 corresponding to a crossing circle gives a tetrahedron when coned to
 the ideal vertex $T^2 \times \{1\}$ on each torihedra.  Since there
 are $c$ crossing circles, this gives $c$ bowties, hence $2c$
 triangular shaded faces and hence $4c$ tetrahedra.  The cones on the
 white faces in each torihedra can be glued to make bipyramids on the
 white faces. These bipyramids can then be stellated into
 tetrahedra. Hence the the number of tetrahedra coming from stellated
 bipyramids equals the number of edges of all the white faces. Since
 an edge of a white face is shared with an edge of a black triangle,
 this equals the number of edges of the torihedral graph, which has $6c$ edges. 
Hence the bipyramids on the
 white faces decompose into $6c$ tetrahedra. Thus total count of
 tetrahedra is $4c+6c=10c$. Since the volume of an ideal tetrahedron is
 bounded by the volume of the regular ideal tetrahedron $\vtet$, the
 upper bound now follows. \qed

\begin{remark}
  In Proposition \ref{prop:falSquareWeave} below we show that our
  upper bound is sharp by showing that the fully augmented square
  weave achieves the upper bound.
\end{remark}


\section{Volume density convergence conjecture}
\label{sec:voldenconvconj}

\subsection{Volume density and its spectrum} In this
section we discuss volume density of fully augmented
links in $\Sp^3$, its spectrum and asymptotic behavior. 
Champanerkar, Kofman and Purcell \cite{CKPgmax} defined volume density
of a hyperbolic link in $\Sp^3$ as the ratio of the volume of the link
complement to its crossing number, and studied the asymptotic behavior
of the volume density for sequences of alternating links which
diagrammatically converge to a biperiodic alternating link.

For a hyperbolic link $L$ in $\Sp^3$, let $\vol(L)$ denote the
hyperbolic volume of ${\rm S}^3 - L$.  In this section we
assume that all links discussed below are hyperbolic.


\begin{define} 
  Let $L$ be a fully augmented link in $\Sp^3$ with or without
  half-twists. The \emph{volume density of $L$} is defined to be the
  ratio of the volume of $L$ and the number of augmentations,
  i.e. $\vol(L)/a(L)$ where $a(L)$ is the number of augmentations of
  the link $L$. We similarly define the volume density of a fully
  augmented link in $T^2 \times I$.
\end{define}

\begin{remark}
  Adams \cite{CA} showed that the volume of an augmented link with a
  half-twist at the crossing circle of the augmentation is equal to
  the volume without a half-twist. However fully augmented links with
  and without half-twists have different crossing numbers. Hence in
  our definition above we divide by the number of augmentations rather
  than the number of crossings.
\end{remark}

\begin{remark}
  For a fully augmented link without half-twist the crossing number of
  the diagram is $4a(L)$. Thus the volume density of such a fully
  augmented link $L$ is related to the volume density of $L$ as
  defined in \cite{CKPgmax} by a factor of $4$. 
 \end{remark}

Throughout this section and the next we consider fully augmented
links without half-twists.

\begin{example}
  The Borromean rings $B$ has $\vol(B)= 2\voct$ and $a(L)=2$, hence the
  volume density $\vol(B)/a(B)= \voct$.
\end{example}

\begin{define}
The volume density spectrum of fully augmented links in $\Sp^3$
is defined as $\mathcal{S}_{aug}= \{\vol(L)/a(L) : L$ fully augmented link in $\Sp^3 \}$.
\end{define}

\begin{prop} \label{thm:falVolumeDensity}
The volume density spectrum $\mathcal{S}_{aug} \subset [\voct, 10\vtet)$.
\end{prop}

{\it Proof.} Let $L$ be a fully augmented link then by Proposition 3.8
of \cite{JessicaP} the volume of $L$ is at least $2\voct
(a(L)-1)$.
Since $L$ is hyperbolic, $a(L) \geq 2$, which implies that the volume
density
$$\frac{\vol(L)}{a(L)} \geq \frac{2\voct \cdot a(L)}{a(L)} -
\frac{2\voct}{a(L)} > 2\voct \bigg(1 - \frac{1}{a(L)}\bigg) \geq
\voct.$$
Since volume density of the Borromean Rings is $\voct$, the lower
bound is realized. Agol and D. Thurston \cite{Lackenby} showed that
$\vol(L) \leq 10\vtet (a(L)-1)$. Hence the the volume density of $L$
is at most $10\vtet$.  \qed


 \begin{figure}
 \centering
 \begin{tabular}{cccc}
 \includegraphics [height=3cm]{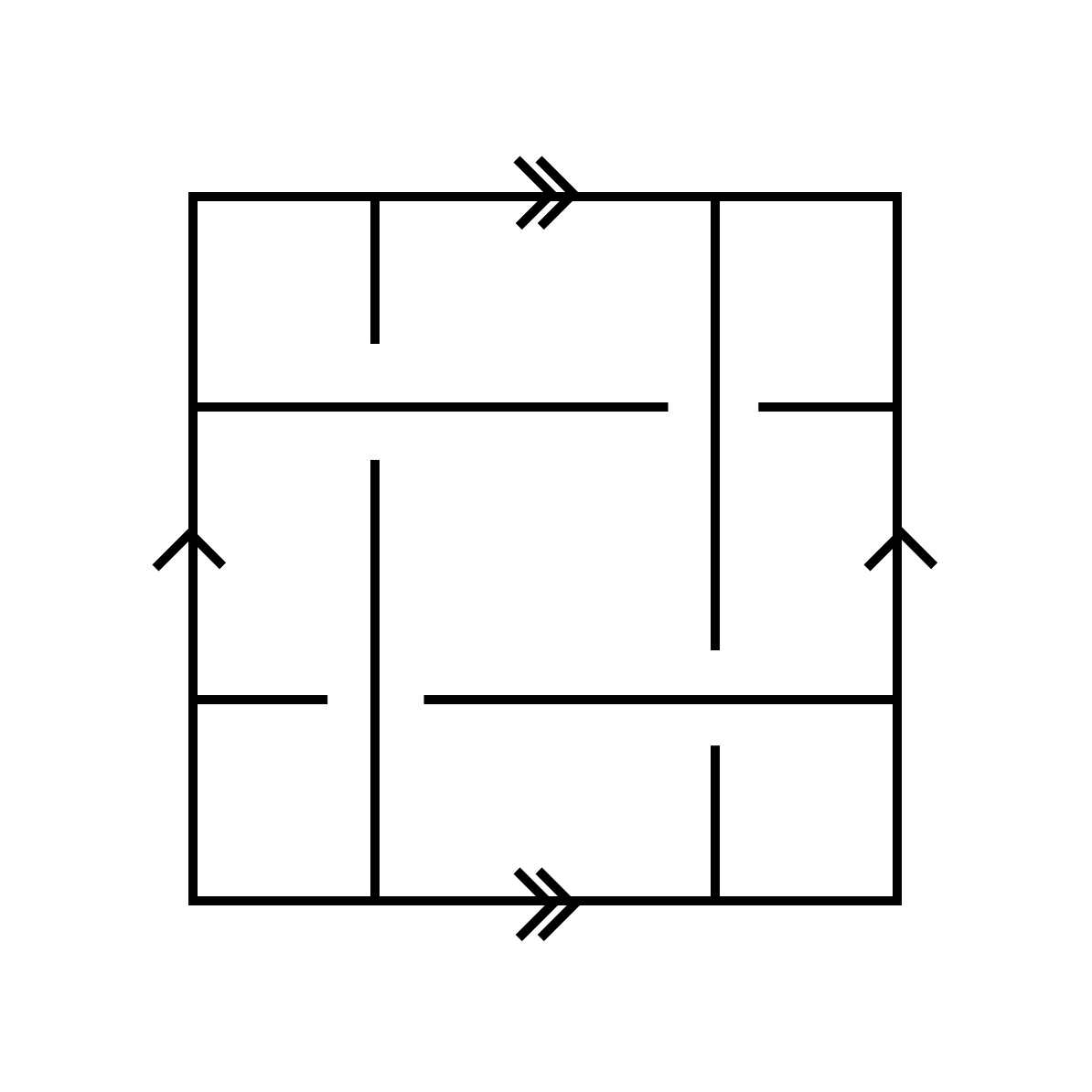}&
 \includegraphics [height=3cm]{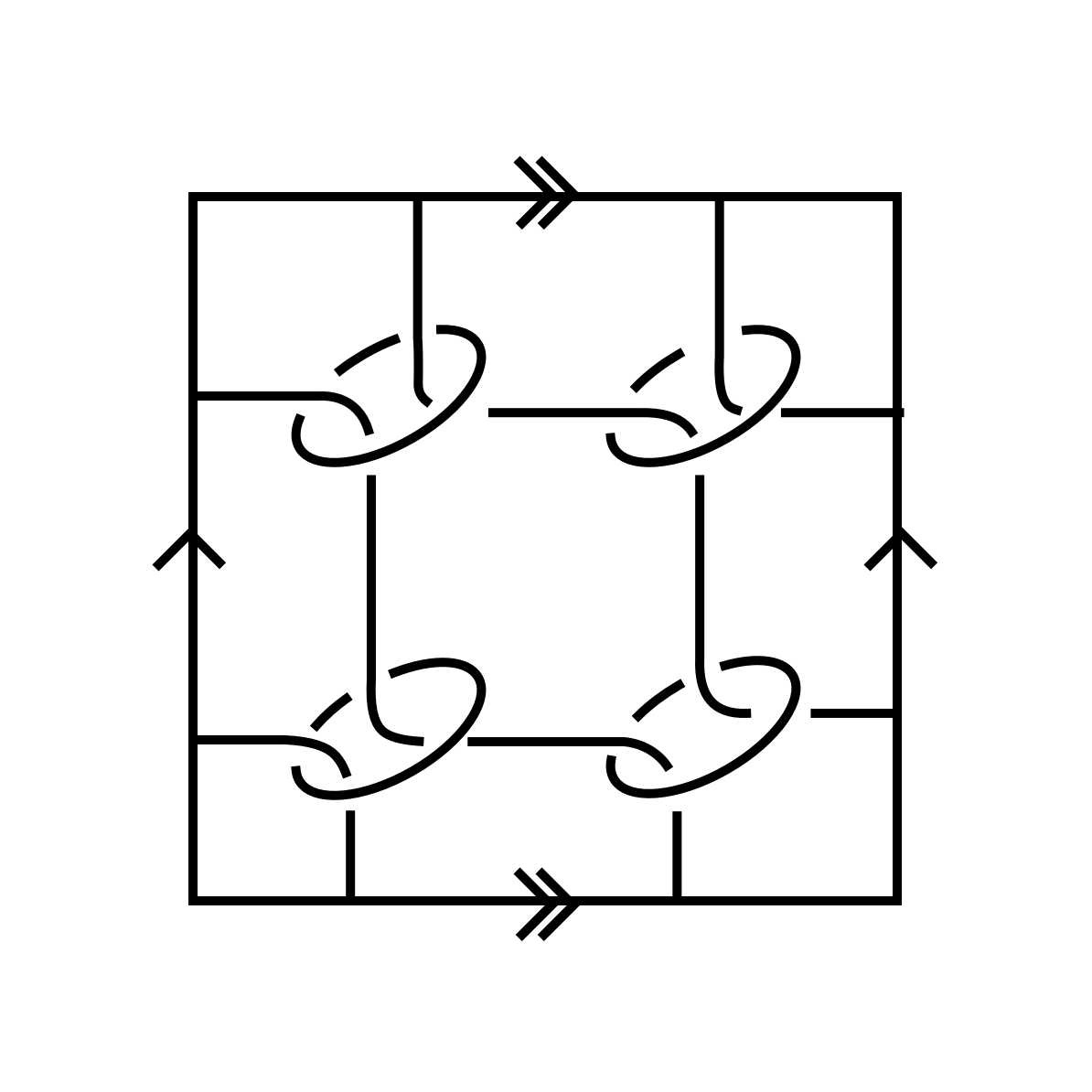}&
 \includegraphics [height=3cm]{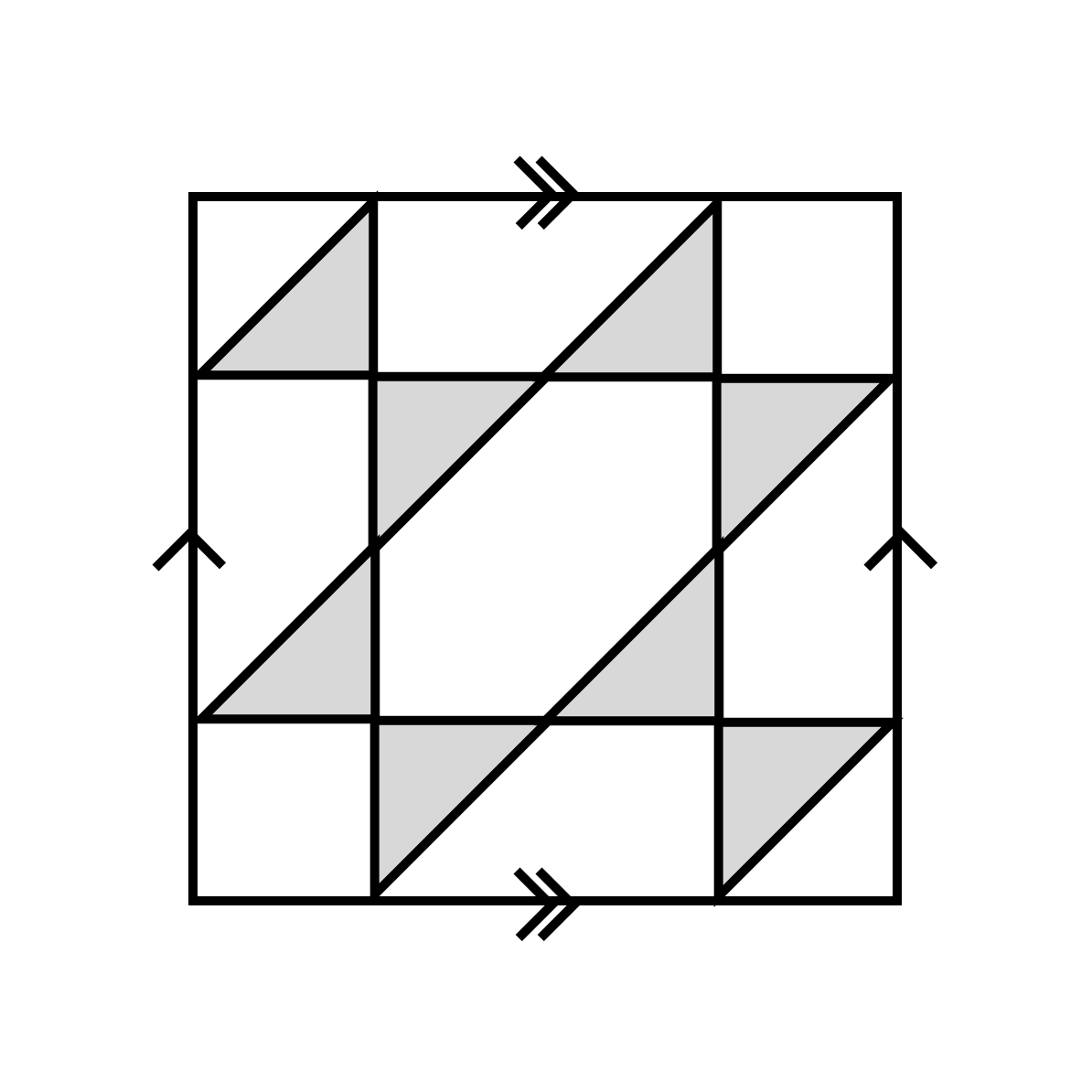}&
  \includegraphics [height=3cm]{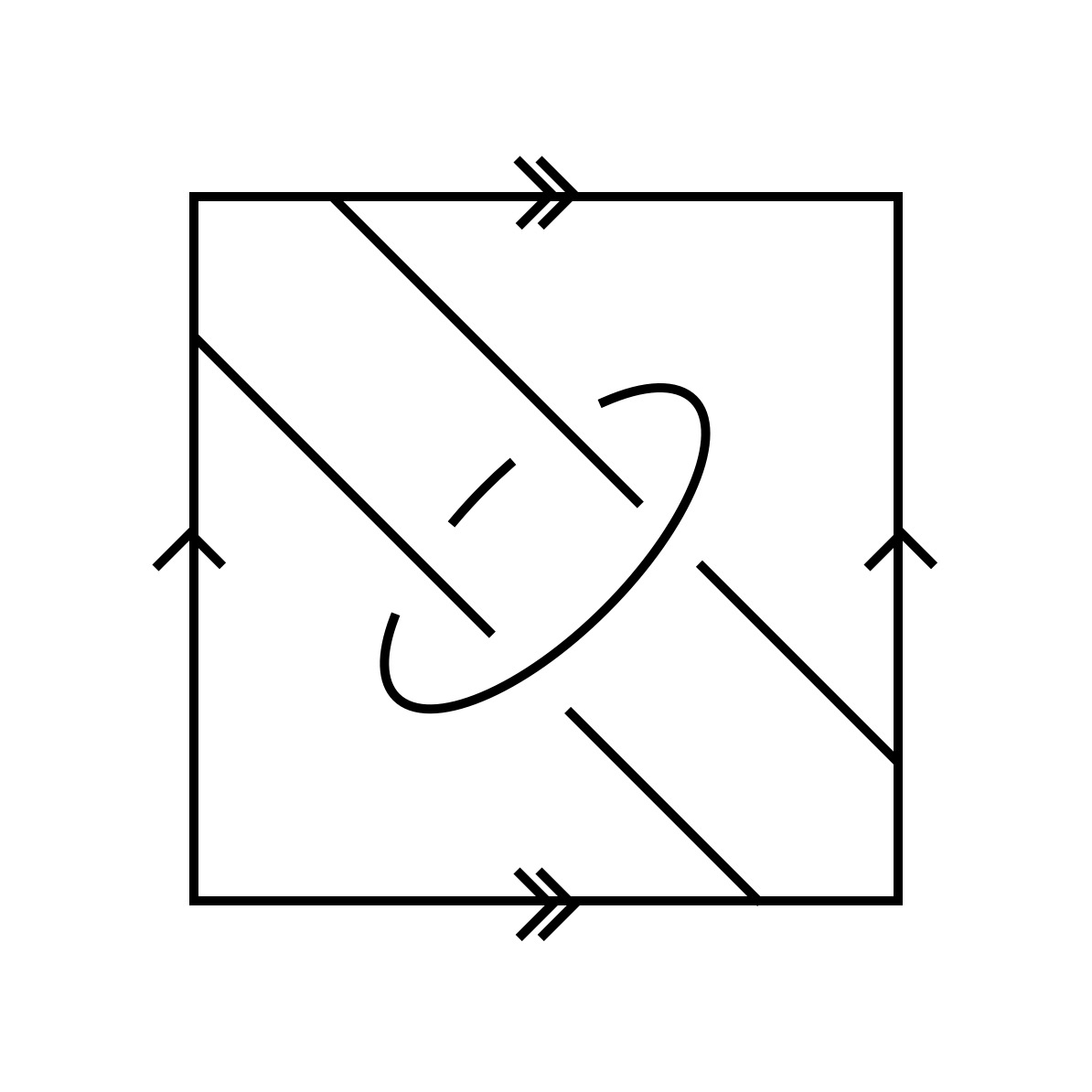}\\
 (a)&(b)&(c)&(d)\\
 \end{tabular}
 \caption{(a) Fundamental domain of the square weave
   $\mathcal{W}$ (b) Fundamental domain of the fully augmented
   square weave, denoted $W_f$. (c) Bow-Tie graph $\Gamma_{W_f}$ of
   the square weave on the left. (d) A quotient of $W_f$ with same
   volume as the triaxial link.}
 \label{fig:falSquareWeave}
 \end{figure}

 We show below that $10\vtet$ occurs as a volume density of the fully
 augmented square weave.  Let $\W_f$ denote the fully augmented square
 weave as in Figure \ref{fig:falSquareWeave}(b).

\begin{prop} \label{prop:falSquareWeave}
$$\frac{\vol(T^2 \times I - W_f)}{ a(W_f) }= 10\vtet $$ 
\end{prop}

{\it Proof.}  A four-fold quotient of $\W_f$ as shown in Figure
\ref{fig:falSquareWeave}(d) was studied in \cite{CFKNPttk}. The
authors proved that the volume of this link complement in the
thickened torus is $10\vtet$. Hence
$\vol(T^2 \times I - W_f) = 40\vtet$, and its volume density 
is $10\vtet$. \qed



\begin{remark}
The quotient of $W_f$ as in Figure \ref{fig:falSquareWeave}(d) has
  the same volume as that of a quotient of a triaxial link which is not a fully augmented link, see Figure
  \ref{fig:fal}(a). However the two links are not the same as they have
  different number of cusps. The triaxial link has 5 cusps, 3 cusps from each link component in the thickened torus and 2 cusps from each link component of the Hopf Link. Whereas, the quotient of $W_f$ in Figure \ref{fig:falSquareWeave}(d) has 4 cusps, 2 from each component of the link in the thickened torus (which includes the crossing circle) and 2 cusps from each link component of the Hopf Link. 


\end{remark}

\subsection{F\o lner convergence}
\label{subsec:Folner}

The volume density of the fully augmented square weave is $10 \vtet$. 
We will prove below that $10\vtet$ is also a limit point of the $\mathcal{S}_{aug}$ by investigating the asymptotic behavior of volume
density of a sequence of fully augmented links in $\Sp^3$ which
diagrammatically converge to the biperiodic fully augmented square weave, as defined below.
We use the the notion of F\o lner convergence, which was
first introduced in \cite{CKPgmax}.
%
%
%
%
%
We begin by modifying the definition of F\o lner convergence given in
\cite{CKPgmax}.

  In \cite{CKPgmax} the authors used the Tait graph (checkerboard graph)
of alternating links to define F\o lner convergence. We will use  bow-tie graphs to define F\o lner convergence 
for fully augmented links, see Definition \ref{def:bowtieGraph} and Proposition 2.2 of \cite{JessicaP}. 



\begin{figure}
\centering
\begin{tabular}{ccc}
\includegraphics [height=3cm]{squareWeave}&
\includegraphics [height=3cm]{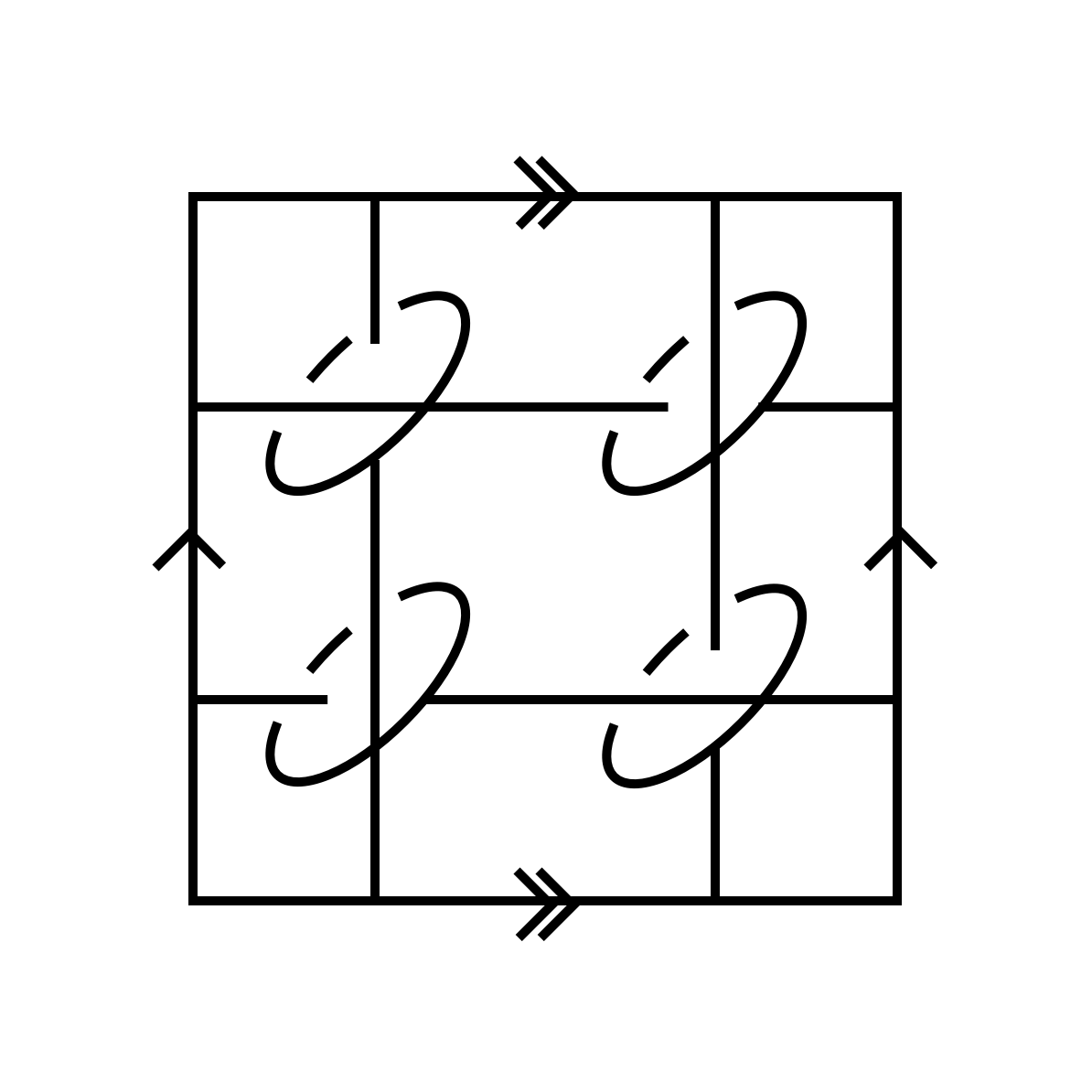}&
\includegraphics [height=3cm]{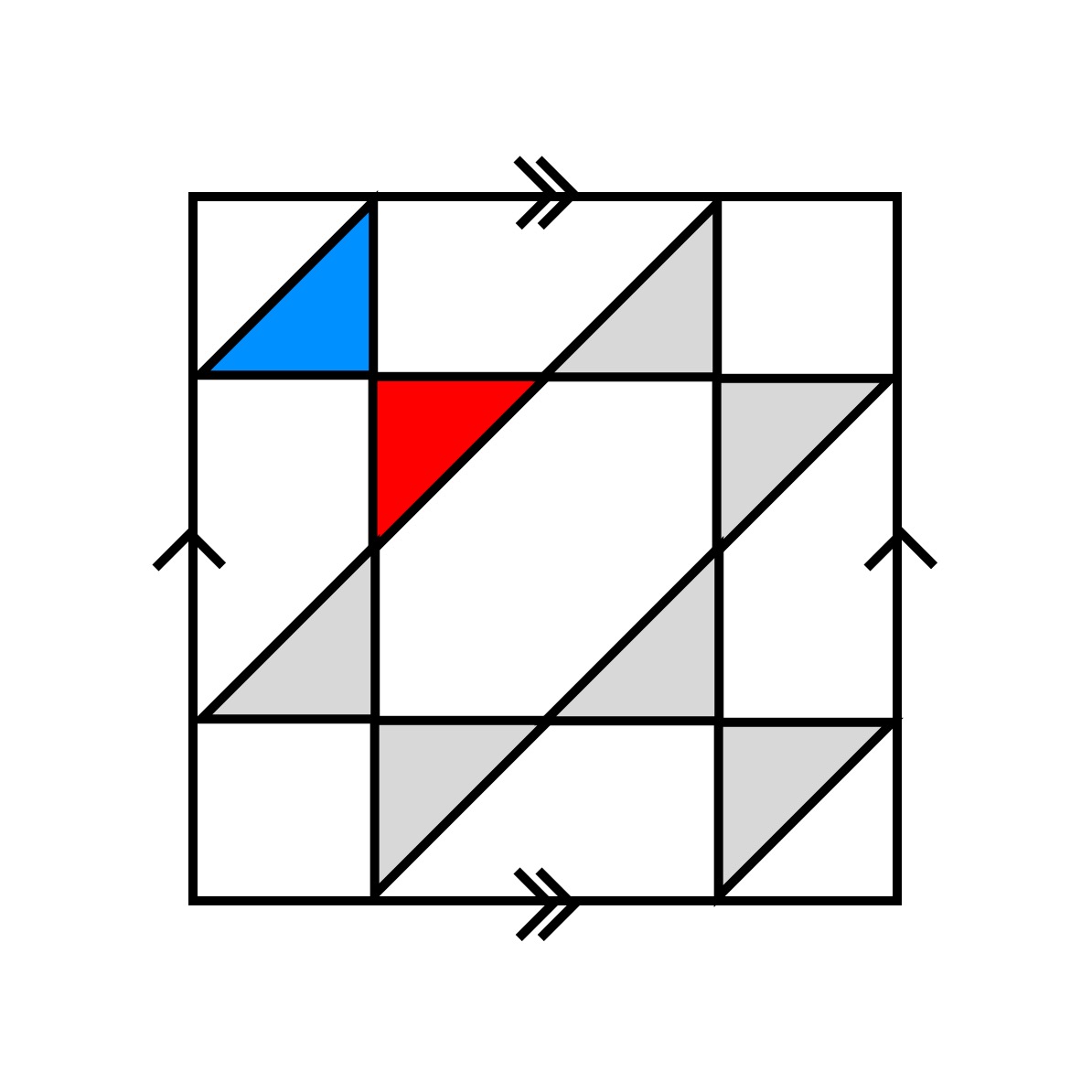}
\includegraphics [height=3cm]{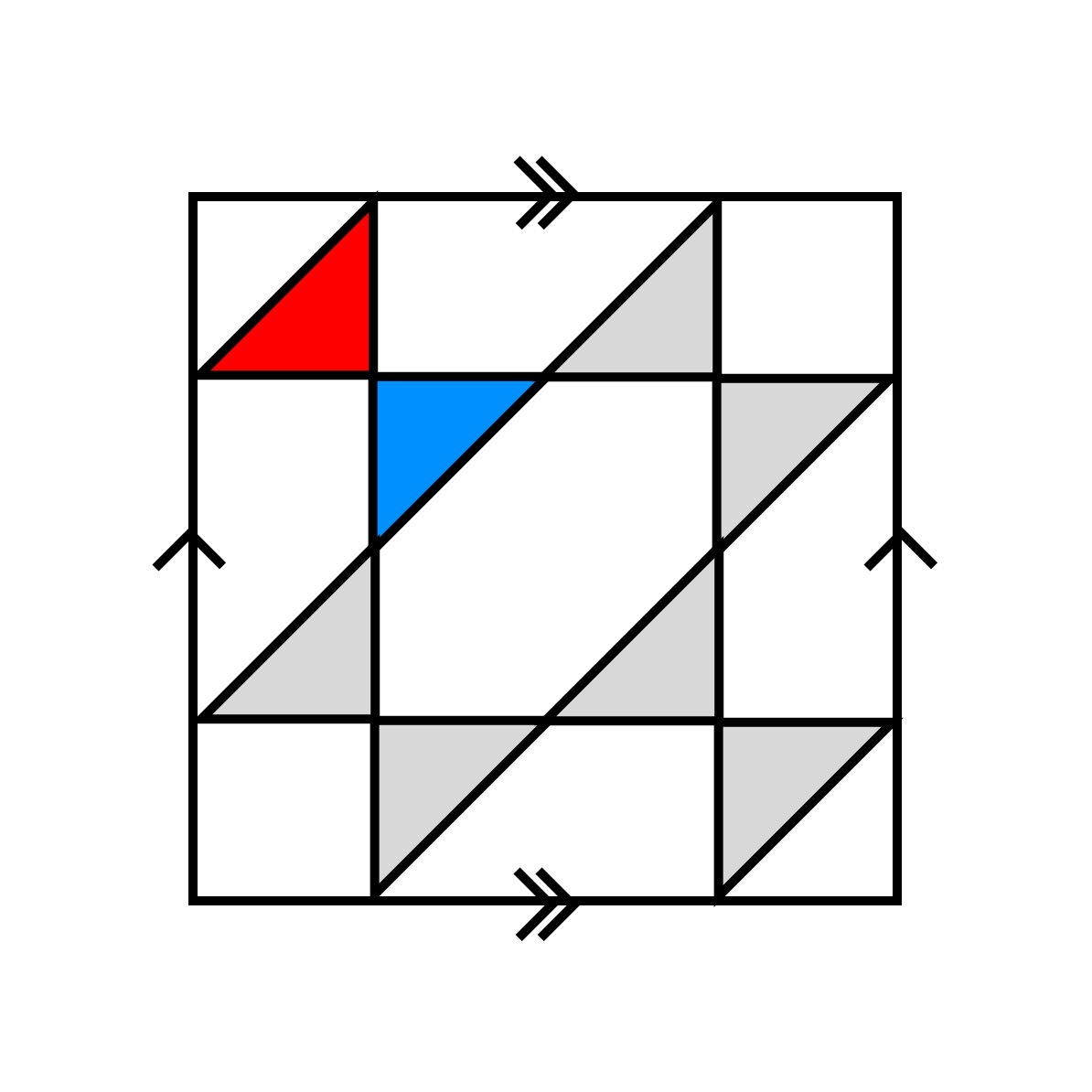}\\
(a)&(b)&(c)
\end{tabular}
\caption{(a) The quotient of the square weave. (b) $W_f$ with half twists at each crossing circle. (c) the bow-tie graph with blue (red) face bow-tie of the top torihedron being glued to a blue (red) face of the bottom torihedron}
\label{fig:halftwist}
\end{figure}

\begin{figure}
\centering
\begin{tabular}{ccc}
\includegraphics [height=3cm]{squareWeave}&
\includegraphics [height=3cm]{WFNOhalftwist}&
\includegraphics [height=3cm]{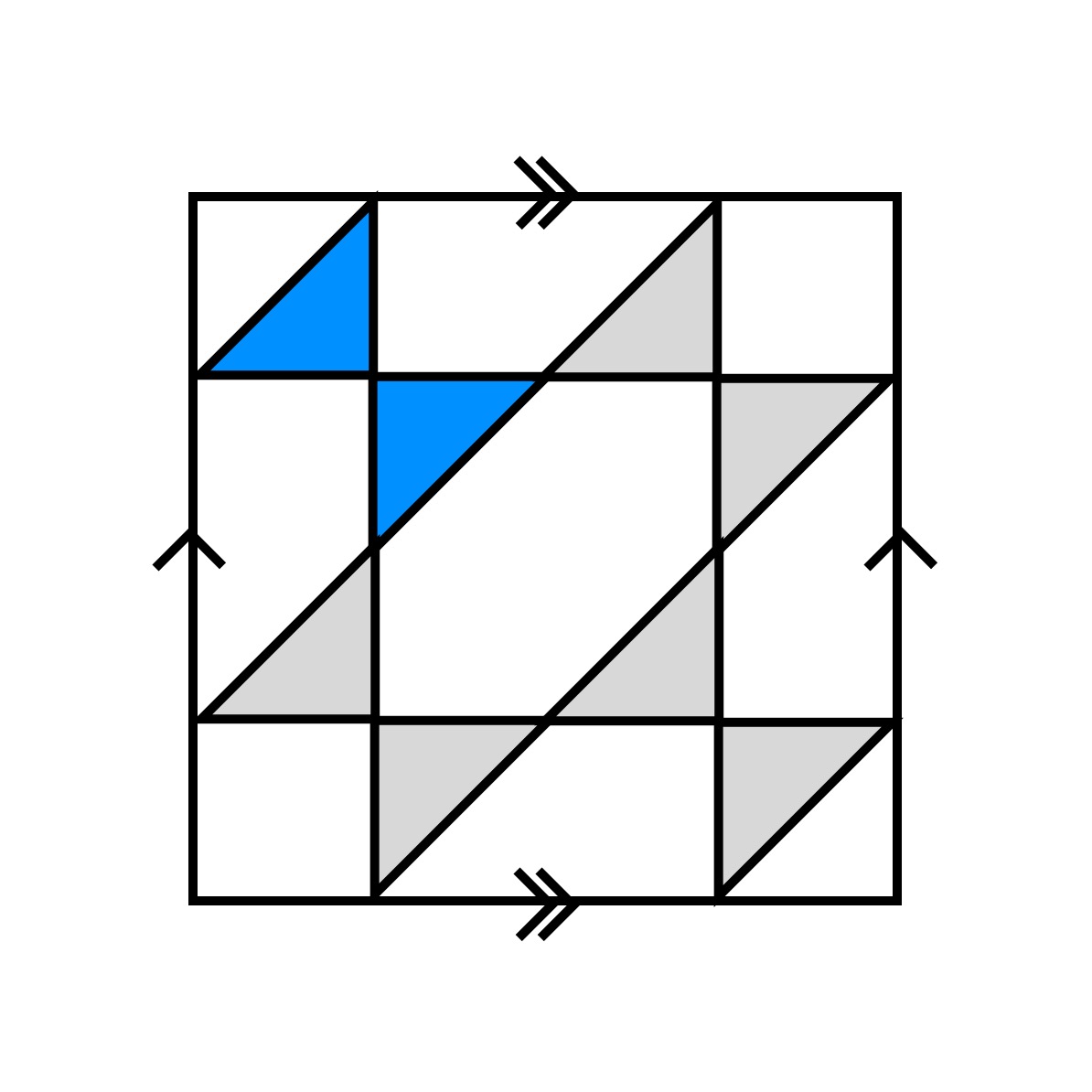}
\includegraphics [height=3cm]{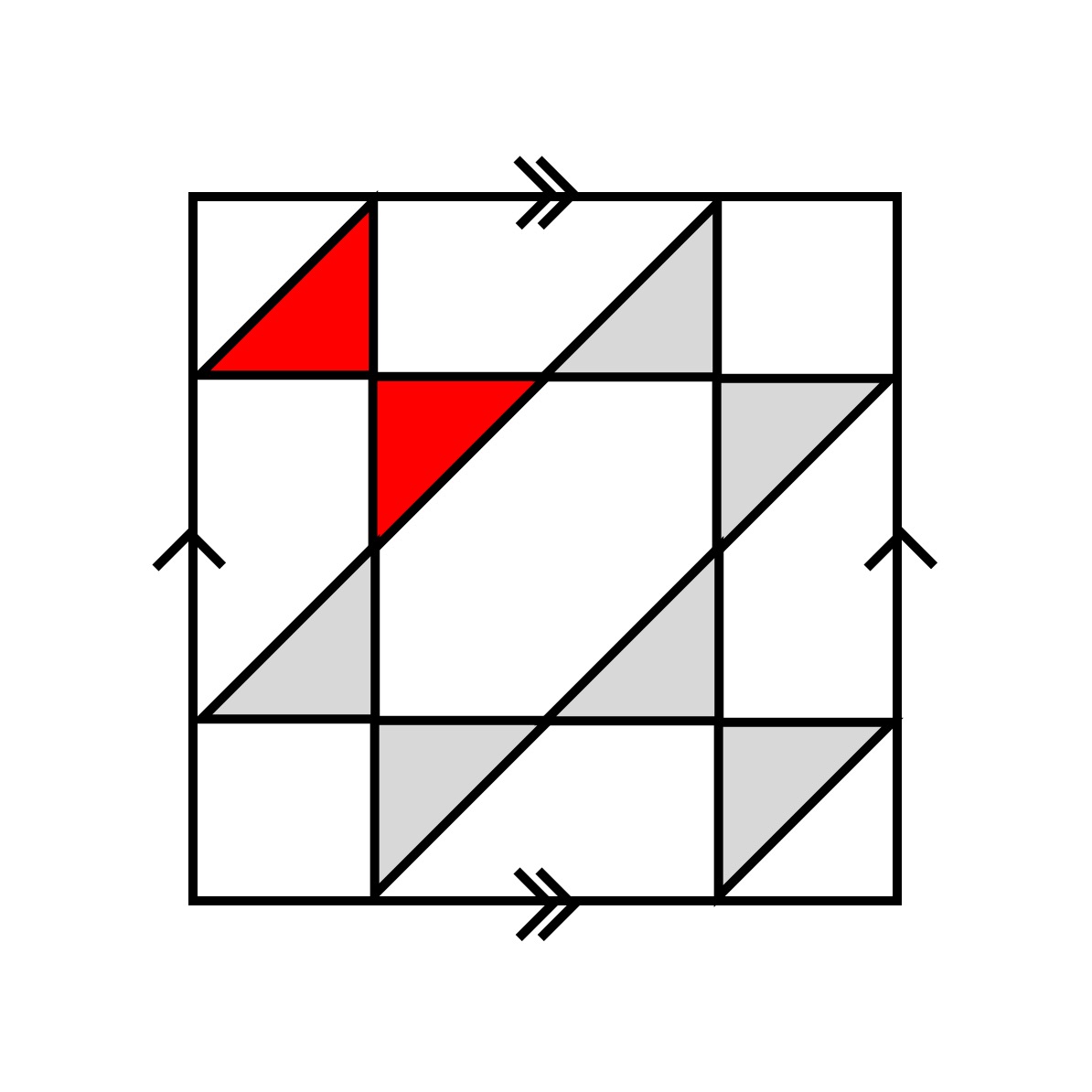}\\
(a)&(b)&(c)
\end{tabular}
\caption{(a) The quotient of the square weave. (b) $W_f$ with no half twists at each crossing circle. (c) the bow-tie graph with blue (red) face bow-tie of the top torihedron being glued to a blue (red) face of the bottom torihedron}
\label{fig:nohalftwist}
\end{figure}

\begin{figure}
\centering
\includegraphics [height=3cm]{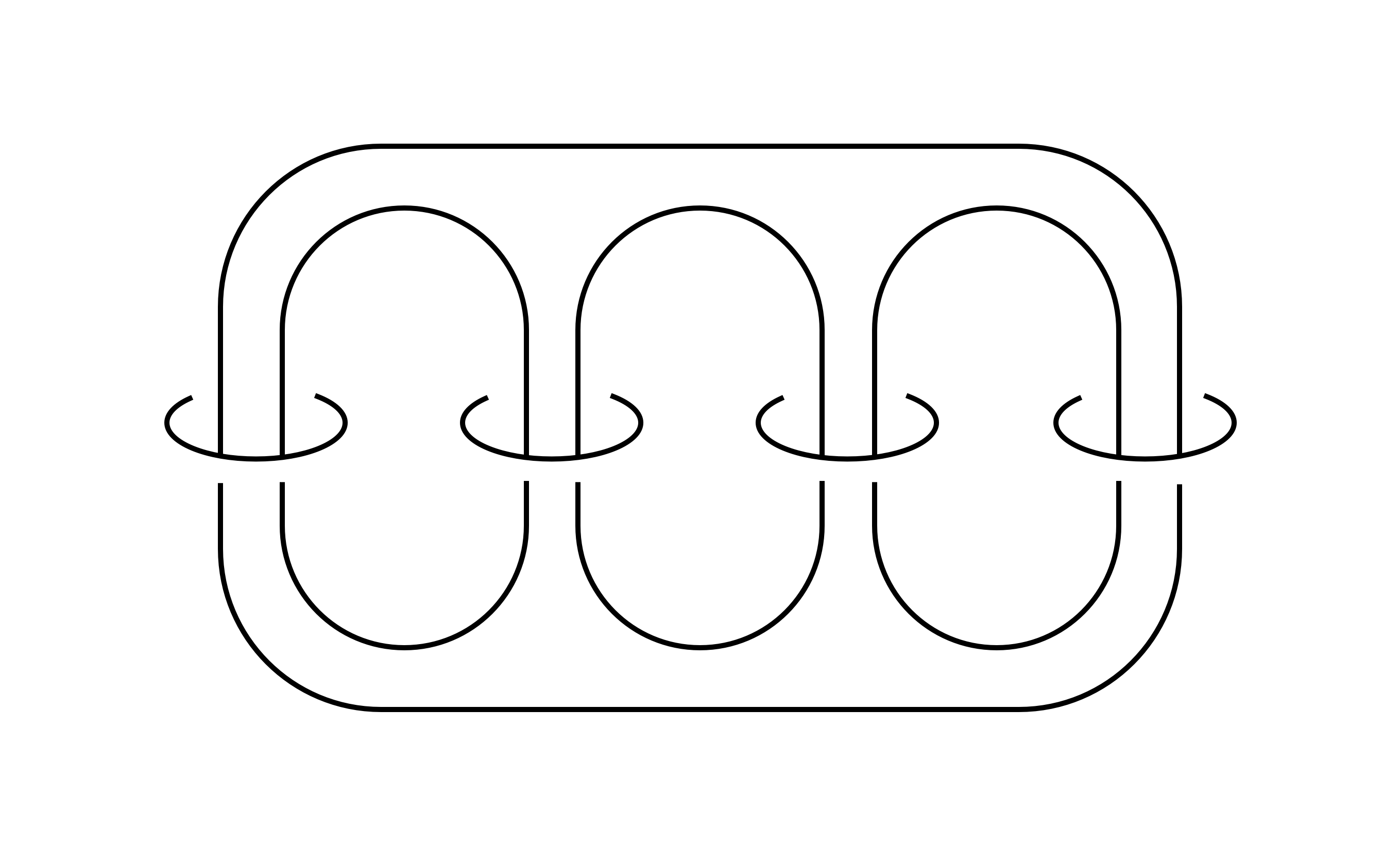}
\includegraphics [height=3cm]{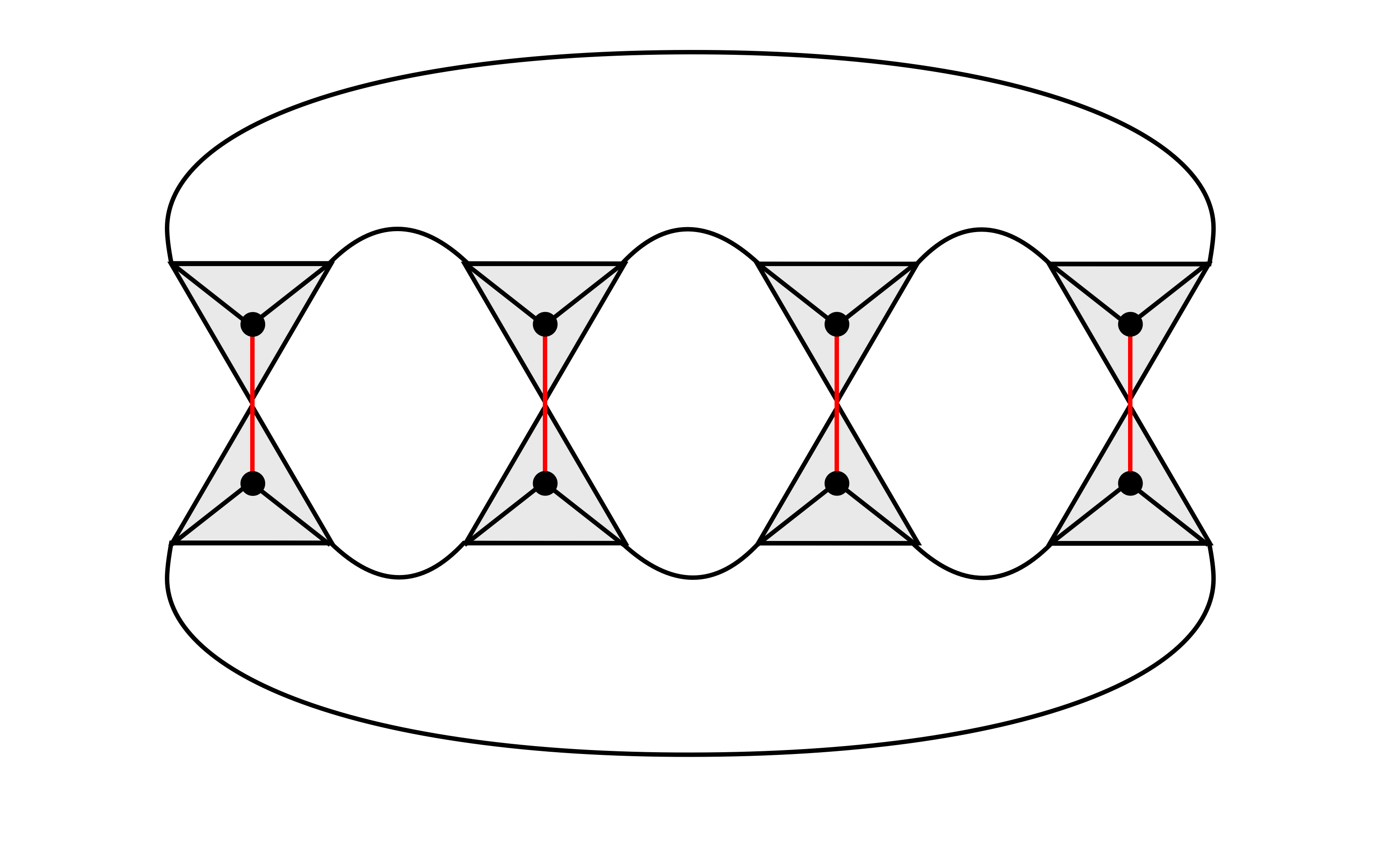}\\
\includegraphics [height=3cm]{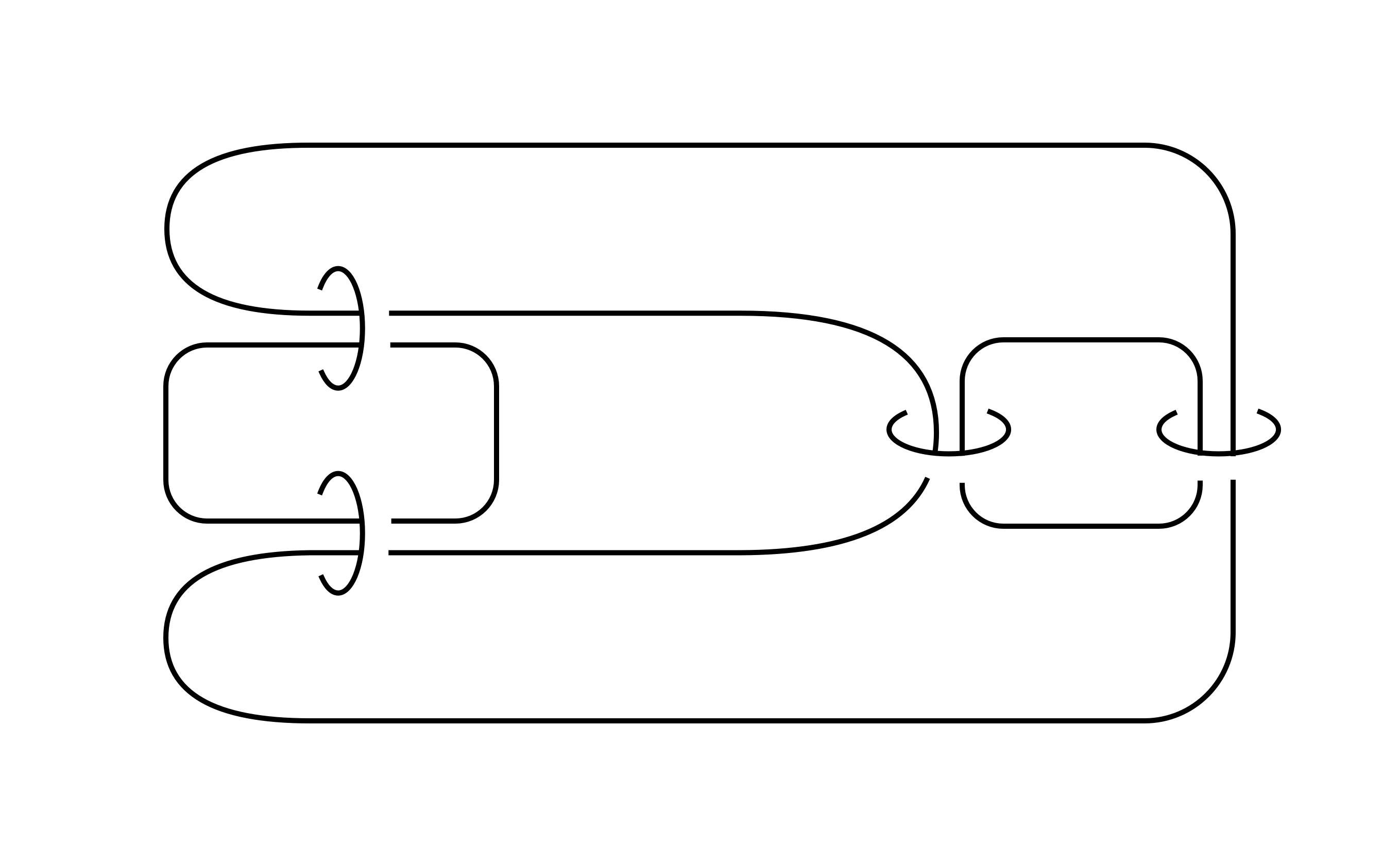}
\includegraphics [height=3cm]{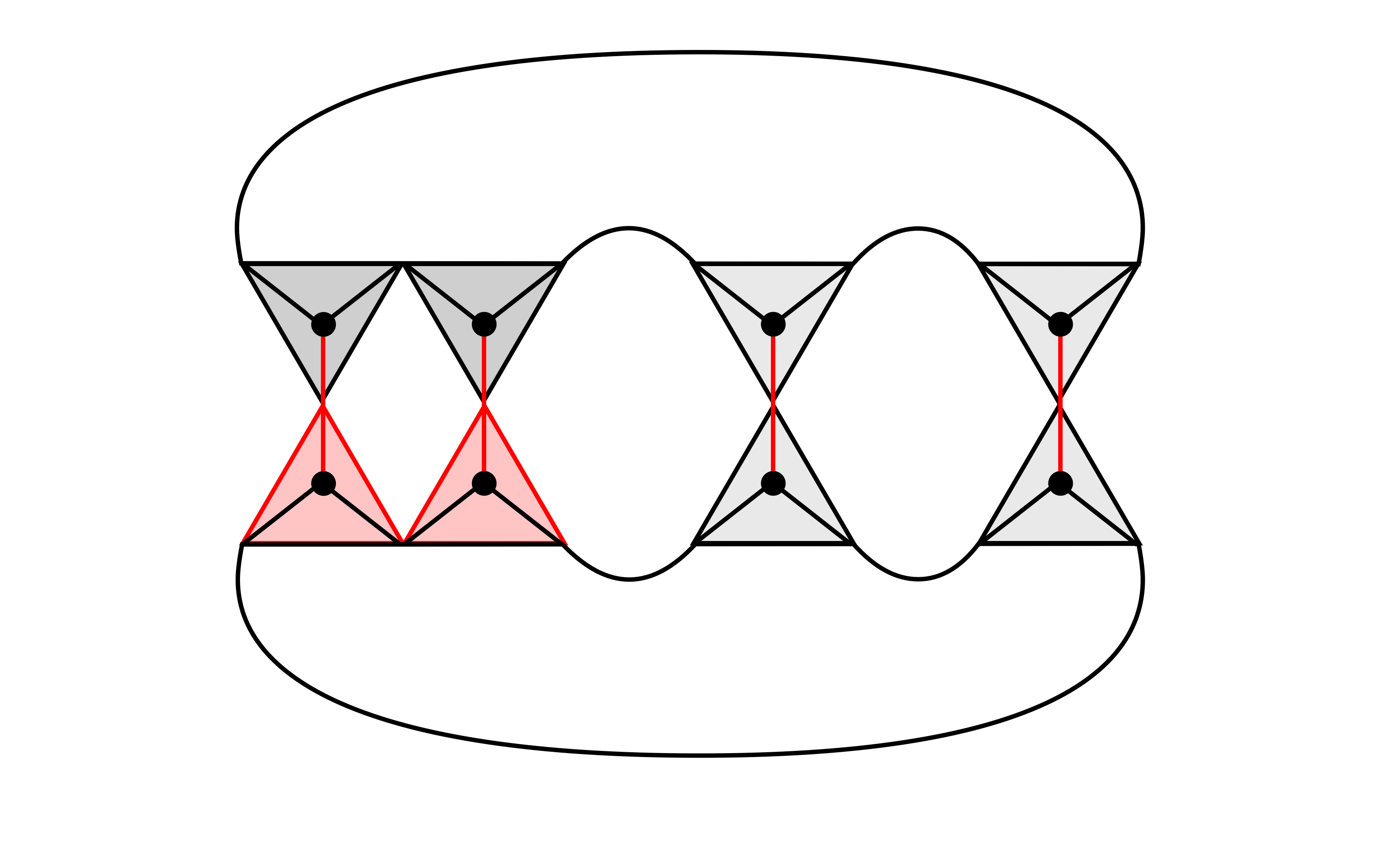}
\includegraphics [height=3cm]{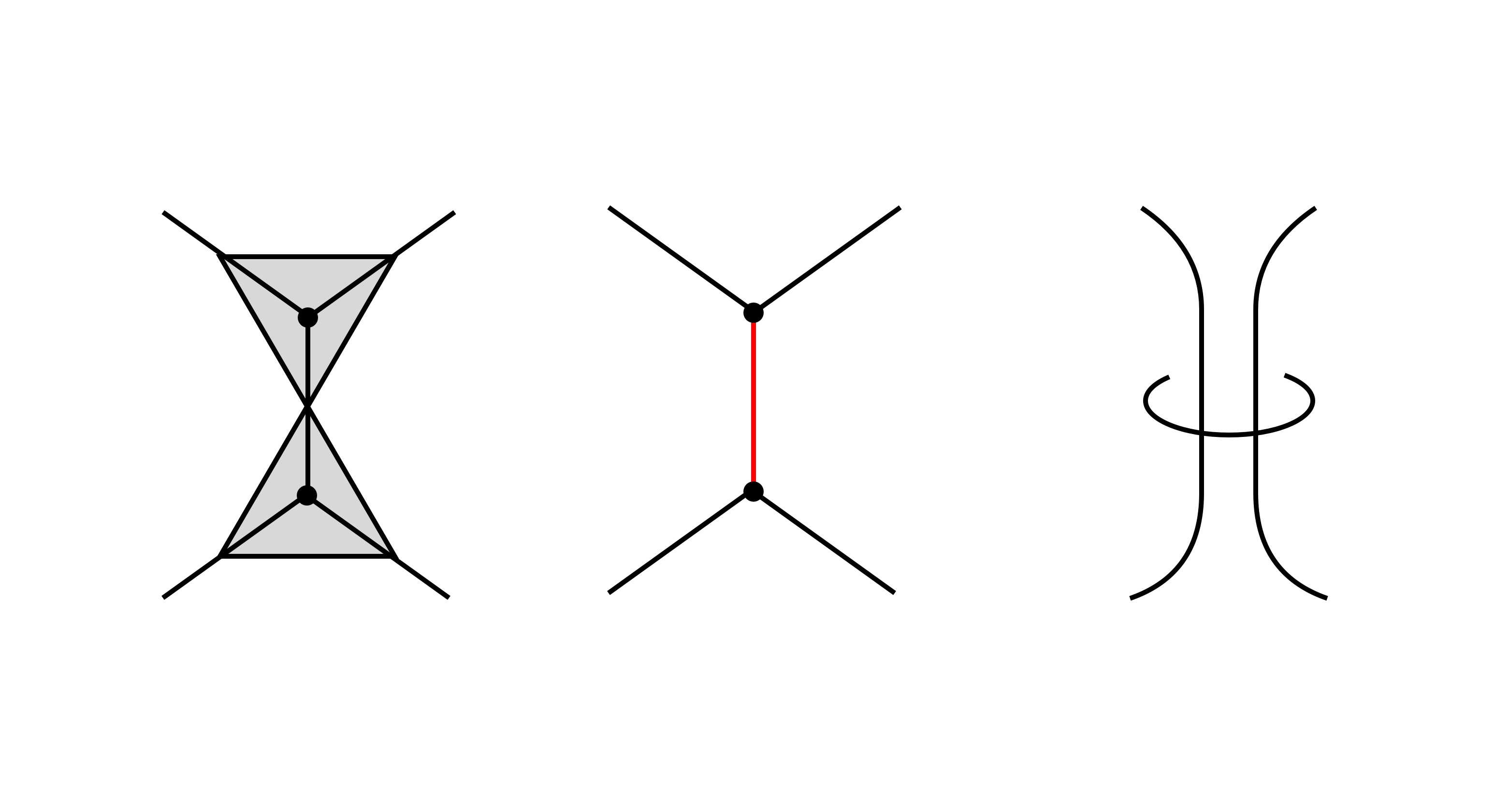}
\caption{The fully augmented links above and below are two different links with the same bow-tie graph with different pairing information.}
\label{fig:diffFAL}
\end{figure}

\begin{define} \label{def:geometricFolner} Let $\mathcal{L}$ be a
  biperiodic fully augmented link. We will say that a sequence of
  fully augmented links $\{K_n\}$ in $\Sp^3$ \emph{F\o lner converges
    almost everywhere geometrically} to $\mathcal{L}$, denoted by
  $K_n \xrightarrow{\text{GF}} \mathcal{L}$, if the respective bow-tie
  graphs $\{\Gamma_{K_n}\}$ and $\Gamma_{\mathcal{L}}$ satisfy the
  following: there are subgraphs $G_n \subset \Gamma_{K_n}$ such that
\begin{enumerate}
\item $G_n \subset G_{n+1}$, and $\cup G_n = \Gamma_{\mathcal{L}}$,
\item $\displaystyle {\lim_{n \to \infty}} |\partial G_n| / |G_n| = 0$, where $| . |$ denotes the number of vertices, and $\partial G_n \subset \Gamma_{\mathcal{L}}$ consists of the vertices of $G_n$ that share an edge in $\Gamma_{\mathcal{L}}$ with a vertex not in $G_n$,
\item $G_n \subset \Gamma_{\mathcal{L}} \cap (n \Lambda)$, where $n \Lambda$ represents $n^2$ copies of the fundamental domain for the lattice $\Lambda$ such that $L = \mathcal{L}/ \Lambda$,
\item $\displaystyle {\lim_{n \to \infty}} |G_n|/ 3a(K_n) = 1$. 
\end{enumerate}
\end{define}

\begin{remark}
  The number 3 appears in the denominator in the last condition for
  the definition of F\o lner convergence because the number of vertices of
  the bow-tie polyhedron for $K_n$ equals three times the number of
  augmentations. To see this note that every bow-tie shares two
  vertices with another bow-tie and hence contributes $3$
  vertices to the graph. Since each bow-tie corresponds to a crossing
  circle, the number of vertices of the graph is $3a(K)$.
\end{remark}

\begin{remark}
Note that many fully augmented links can have the same bow-tie
graph. For example, a fully augmented link with and without
 half-twists have the same bow-tie graph but different gluing. See Figure
\ref{fig:halftwist} and Figure \ref{fig:nohalftwist}. 
%
Another example of this is when the bow-tie graphs are same but with
different pairing of triangles. See Figure \ref{fig:diffFAL} for an
example of two links with same bow-tie graphs but different pairings.
In our definition above, we are using only  the polyhedral graphs 
but not the pairing information of the bow-ties. Hence 
we call our F\o lner convergence geometric. This has the advantage of having 
many more sequences converging to a given biperiodic fully augmented link. 
\end{remark}

\subsection{Volume Density Conjecture}\label{sec:VolDetConj}

\begin{conjecture} (Volume Density Conjecture) \label{con:volumeDensityConjecture}
Let $\mathcal{L}$ be any biperiodic alternating link with alternating quotient link $L$. Let $\{K_n\}$ be a sequence of alternating hyperbolic links such that $K_n$ F\o lner converges to $\mathcal{L}$. Then $$\displaystyle {\lim_{n \to \infty}}  \frac{ \vol(K_n) }{c(K_n)} = \frac{\vol((T^2 \times I) - L)}{c(L)}.$$ 
\end{conjecture}

Champanerkar, Kofman and Purcell proved this conjecture when
$\mathcal{L}$ is the square weave \cite{CKPgmax} and the triaxial link
\cite{CKPtriax}. 
by finding upper and lower bounds on $\vol(K_n)$
such that for a sequence of alternating links
$K_n \xrightarrow{\text{F}} \mathcal{L}$, these bounds
are equal in the limit. One of the key tools in their proof is the use
of right-angled circle patterns. Using the right-angled
decomposition of fully augmented link complements in $\Sp^3$ we
construct  right-angled circle patterns, and use these
to prove the Volume
Density Conjecture for fully augmented links in $\Sp^3$.

The idea is as follows: As described in \cite{JessicaP} each
hyperbolic fully augmented link complement in $\Sp^3$ can be
decomposed into two right-angled ideal polyhedra which are described
by a right-angled circle pattern. By Theorem \ref{thm:fal} each
torihedra of the bow-tie torihedral decomposition are right-angled
described by another right-angled circle pattern on the torus. The
$\Z \times \Z$ lift of this circle pattern is the circle pattern
associated to $\mathcal{L}$. 
We show below that when a sequence of fully augmented links $K_n$ converges 
to $\mathcal{L}$, 
 $K_n \xrightarrow{\text{GF}} \mathcal{L}$, the circle
pattern for $K_n$ converges to infinite circle pattern for
$\mathcal{L}$. As a consequence we obtain the volume density
convergence.

In order to work with circle patterns and convergence of circle
patterns, we recall the following definitions from \cite{atkinson}:
\begin{define} 
A \emph{disk pattern} is a collection of closed round disks in the plane such that no disk is the Hausdorff limit of a sequence of distinct disks and such that the boundary of any disk is not contained in the union of two other disks.
\end{define}
\begin{define} 
  A \emph{simply connected} disk pattern is a disk pattern in the
  plane so that the union of the disks is simply connected.
\end{define}

\indent Let $D$ be a disk pattern in $\C$. Let $G(D)$ be the graph
with a vertex for each disk and an edge between any two vertices when
the corresponding disks overlap. The graph $G(D)$ inherits an
embedding in the plane from the disk pattern and we will identify
$G(D)$ with its plane embedding. A face of $G(D)$ is an unbounded
component of the complement of $G(D)$ in the plane. We can label the edges of
$G(D)$ with the angles between the intersecting disks.

\begin{define}
  A disk pattern $D$ is called an \emph{ideal disk pattern} if the
  labels of edges of $G(D)$ are in the interval $(0,\pi/2]$ and the
  labels around each triangle or quadrilateral in $G(D)$ sum to $\pi$
  or $2\pi$ respectively.
\end{define}

It is clear that ideal disk patterns in $\C$ correspond to ideal
polyhedra in $\mathbb{H}^3$, with the disks corresponding to the faces of 
the ideal polyhedron. 

\begin{define} \label{def:generation}
Let $D$ and $D'$ be disk patterns. Give
  $G(D)$ and $G(D')$ the path metric in which each edge has length
  1. For disks $d$ in $D$ and $d'$ in $D'$, we say $(D,d)$ and
  $(D',d')$ \emph{agree to generation n} if the balls of radius $n$
  centered at vertices corresponding to $d$ and $d'$ admit a graph
  isomorphism, with labels on edges preserved.
 \end{define}

\begin{figure}
 \centering
 \includegraphics [height=6cm]{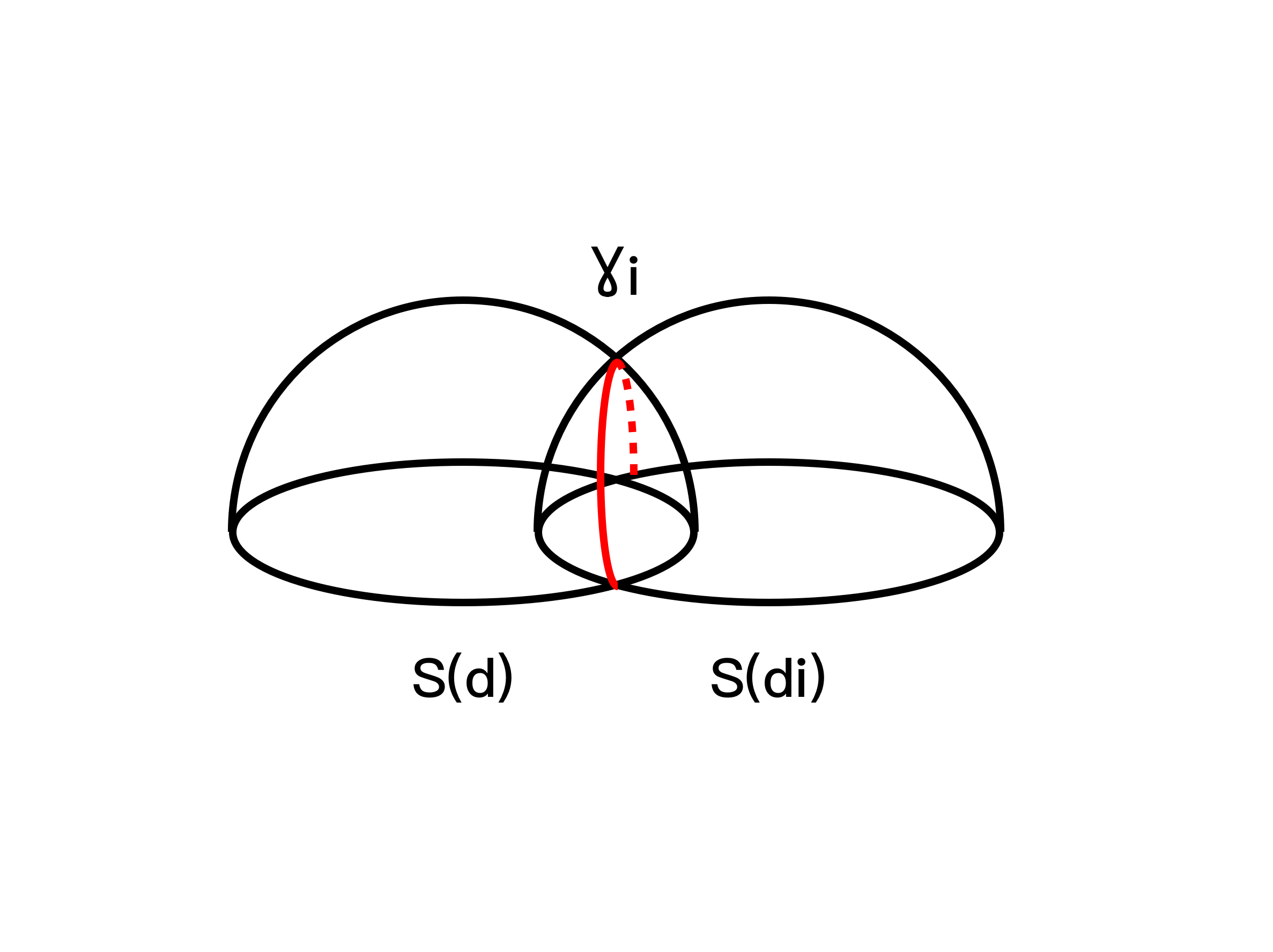}
\includegraphics [height=6cm]{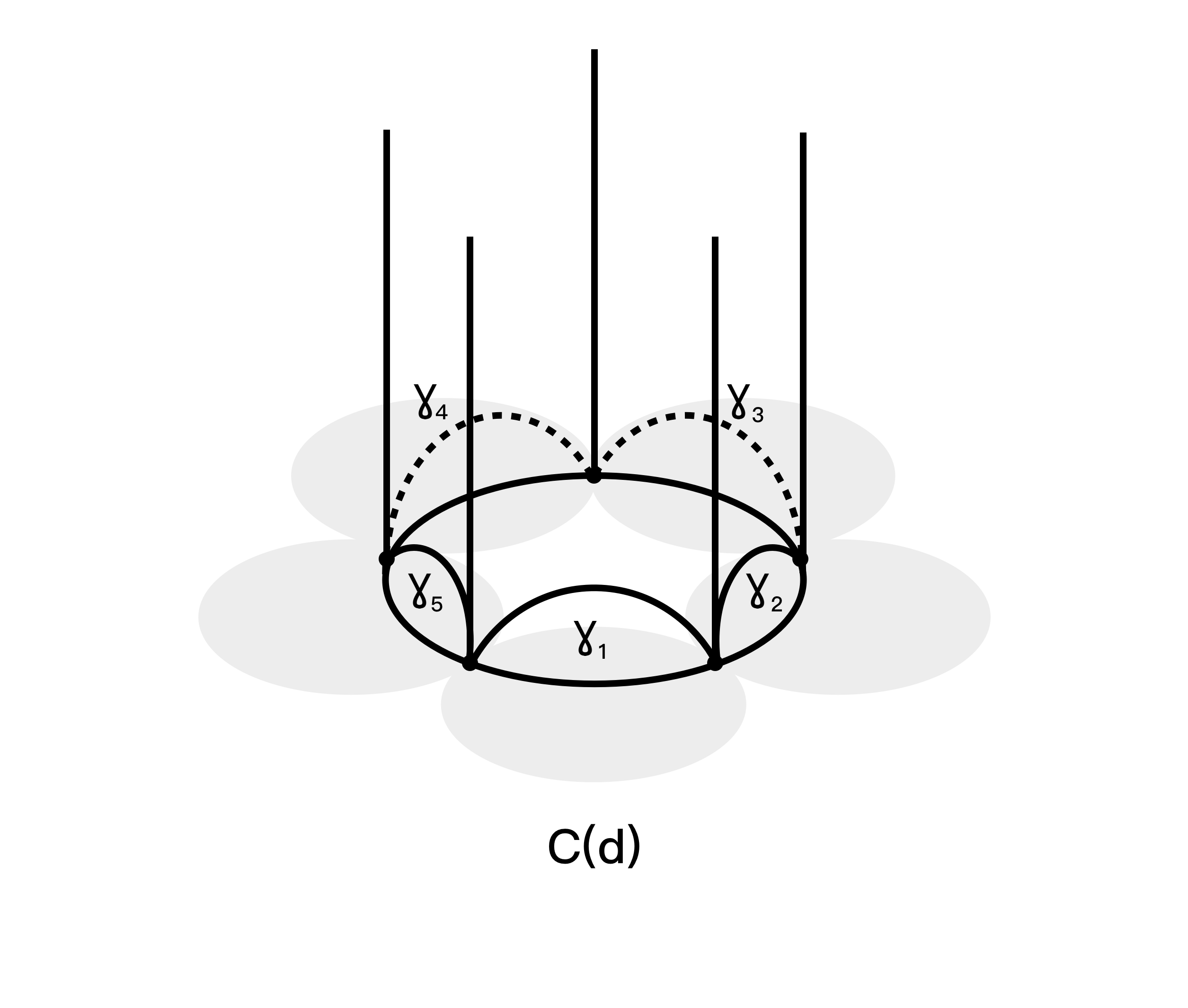}
 \caption{Left: $S(d) \cap S(d')$ Right: C(d)}
 \label{fig:sdCd}
 \end{figure}

 \begin{define} 
   For a disk $d$ in a disk pattern $D$, we let $S(d)$ be the geodesic
   hyperplane in $\HH^3$ whose boundary agrees with that of $d$. That
   is, $S(d)$ is the Euclidean hemisphere in $\HH^3$ with boundary
   coinciding with the boundary of $d$. For a disk pattern coming from
   a right-angled ideal polyhedron, the planes $S(d)$ form the
   boundary faces of the polyhedron. In this case, the disk pattern
   $D$ is \emph{simply connected} and \emph{ideal}, since it
   corresponds to an ideal polyhedron.

   \ \\
   \indent Similarly, for a disk $d$ in $D$, with intersecting
   neighboring disks $d_1,...,d_m$, the intersection
   $S(d) \cap S(d_i)$ is a geodesic $\gamma_i$ in $\HH^3$. The
   geodesics $\gamma_i$ for $i=1, \ldots m$ on $S(d)$ bound an ideal
   polygon in $\HH^3$. The cone of this polygon to the point at infinity 
  is denoted by $C(d)$. See Figure \ref{fig:sdCd}.
 \end{define}

 \begin{define} \label{def:rigid}
 A disk pattern $D$ is said to be \emph{rigid} if $G(D)$ has only triangular and quadrilateral faces and each quadrilateral face has the property that the four corresponding disks of the disk pattern intersect in exactly one point.
 \end{define}
 
 We will use the following lemma proved by Atkinson in
 \cite{atkinson}. 
 
 \begin{lemma} \cite{atkinson} \label{lem:atkinson} Let $D_{\infty}$
   be an infinite rigid disk pattern.  Then there exists a bounded sequence
   $0 \leq \epsilon_l \leq b < \infty$ converging to zero so that if
   $D$ is a simply connected, ideal, rigid finite disk pattern
   containing a disk $d$ so that $(D_{\infty}, d_{\infty})$ and
   $(D,d)$ agree to generation $l$ then
   $$|\vol(C(d)) - \vol(C(d_{\infty}))| \leq \epsilon_l.$$
\end{lemma}

Note that the sequence $\{ \epsilon_l \}$ in above Lemma only depends
on $D_{\infty}$.



We are now ready to prove: 
\begin{theorem} (Volume Density Conjecture for fully augmented
  links) \label{thm:volDetConjFAL}
Let $\mathcal{L}$ be a biperiodic fully augmented link with quotient link $L$. Let $\{K_n\}$ be a sequence of hyperbolic fully augmented links in $\Sp^3$.
Then
$$ K_n \xrightarrow{\text{GF}}  \mathcal{L} \implies \lim_{n \to \infty}  \frac{ \vol(K_n) }{a(K_n)} = \frac{\vol((T^2 \times I) - L)}{a(L)}.$$
\end{theorem}

{\it Proof:} 
Let $P_L$ be the bow-tie torihedron with bow-tie graph $\Gamma_L$ of
$L$. Let $P_{\infty}$ be the infinite polyhedron in $\mathbb{H}^3$
which is the biperiodic lift of $P_L$ with its cone vertex taken to be
$\infty$.  $P_{\infty}$ 
can be seen to be made up of $\Z^2$
copies of an embedding of $P_L$ in $\mathbb{H}^3$ with its cone
vertex taken to be $\infty$, glued according to the biperiodic lift.
Note that since the graph of $P_L$ is the bow-tie graph $\Gamma_L$ of
$L$, which is toroidal, the graph of $P_{\infty}$ is a biperiodic lift
of $\Gamma_L$ and is isomorphic to the bow-tie graph
$\Gamma_{\mathcal{L}}$ coming from $\mathcal{L}$. 
Let $D_{\infty}$ be
the infinite disk pattern coming from the infinite polyhedron
$P_{\infty}$. Since $P_L$ is right-angled torihedron, $P_{\infty}$ is
also right-angled, and hence $D_{\infty}$ is a right-angled disk pattern.

Since $\{K_n\}$ is a sequence of fully augmented links, each $K_n$ is a fully augmented hyperbolic link in $\Sp^3$.
The bow-tie
polyhedron of $K_n$ is a right-angled ideal hyperbolic polyhedron with the same
graph as the bow-tie graph $\Gamma_{K_n}$. The assumption that the
sequence $\{K_n\}$ F\o lner converges almost everywhere geometrically
to $\mathcal{L}$ implies that there are subgraphs
$G_n \subset \Gamma_{K_n}$ which satisfy the conditions of F\o lner
convergence in Definition \ref{def:geometricFolner}.  Hence we can
embed bow-tie polyhedra of $K_n$ in $\mathbb{H}^3$ such that they 
satisfy the following two conditions: (1) Pick a vertex in
$\Gamma_{K_n} - G_n$ to send to infinity, and (2) $G_n \subset G_{n+1}$. 
We denote this polyhedron in $\mathbb{H}^3$ by $P_n$.  First note that $\vol(K_n) = 2\vol(P_n)$. 
Let $v(P_n)$ denote the number of vertices of $P_n$. Since $P_n$ is a 4-valent checkerboard graph whose shaded faces are triangles coming from the bow-ties, one for each augmentation,  every vertex is shared by two triangles. Hence $v(P_n)= 3\cdot 2a(K_n)/2=3a(K_n)$. Therefore, 
 $$ \frac{\vol(K_n)}{3 a(K_n)} = 2 \frac{\vol(P_n)}{v(P_n)}.$$

\begin{figure}
 \centering
 \includegraphics [height=5cm]{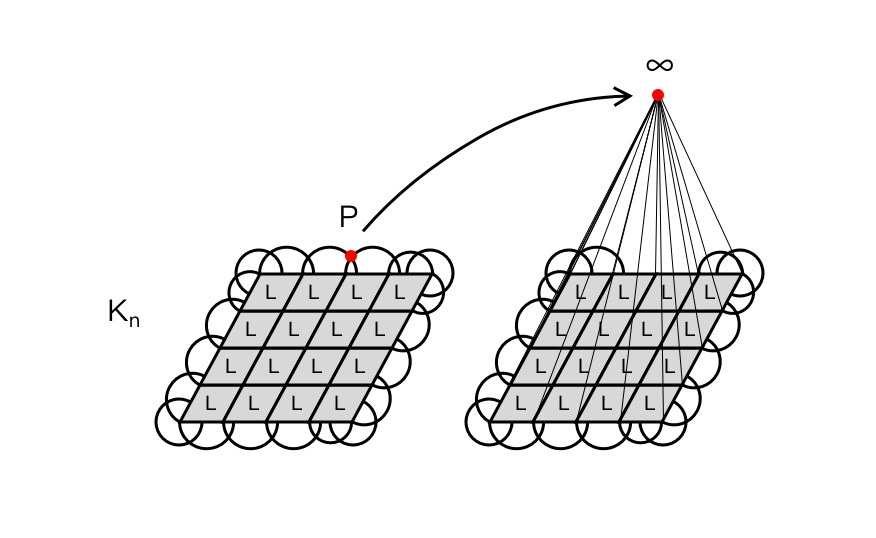}
 \caption{Left: An $n \times n$ copies of fundamental domain of $\Lambda$  
 with an arbitrary closure, and a marked point $P$ on the crossing of the closure. 
 Right: Point $P$ moved to cone point at $\infty$ .}
 \label{fig:nbynL}
 \end{figure}

Let $D_n$ be
the disk pattern of the polyhedron $P_n$. It follows that $D_n$ is a
right-angled, simply connected disk pattern. Since $D_n$ corresponds
to a disk pattern arising from a fully augmented link, 
$D_n$ is rigid (see Definition \ref{def:rigid}, Figure
\ref{fig:nbynL}).
We will now use F\o lner convergence to relate $D_n$ and
$D_{\infty}$.

Let $F_l^n$ be
the set of disks $d$ in $D_n$ so that $(D_n, d)$ agrees to generation
$l$ but not to generation $l+1$ with $(D_{\infty}, d_{\infty})$. For
every positive integer $k$, let $|f_k^n|$ denote the number of faces
of $P_n$ with $k$ sides that are not contained in $\cup_{l} F_l^n$ and
do not meet the point at infinity. By counting vertices we obtain 
$$\sum_k k|f_k^n| \leq 4|\Gamma_{K_n}-G_n|.$$ 

The term $|\Gamma_{K_n}-G_n|$ counts the number of vertices that are in $\Gamma_{K_n}$ but not in $G_n$. Since all the vertices of the graph $\Gamma_{K_n}$ are four valent we get a factor of $4$.  Hence 
$|\Gamma_{K_n}| = v(P_n) = 3a(K_n)$, and  

\begin{equation} \label{eqn:one}
\displaystyle {\lim_{n \to \infty}} \frac{|G_n|}{3a(K_n)} = 1 \implies 
\displaystyle {\lim_{n \to \infty}} \frac{4|\Gamma_{K_n} - G_n|}{v(P_n)}=0 \implies 
\displaystyle {\lim_{n \to \infty}} \frac{\sum_k k|f_k^n|}{v(P_n)}=0.
\end{equation}
 
 \begin{figure}
 \centering
 \includegraphics [height=8cm]{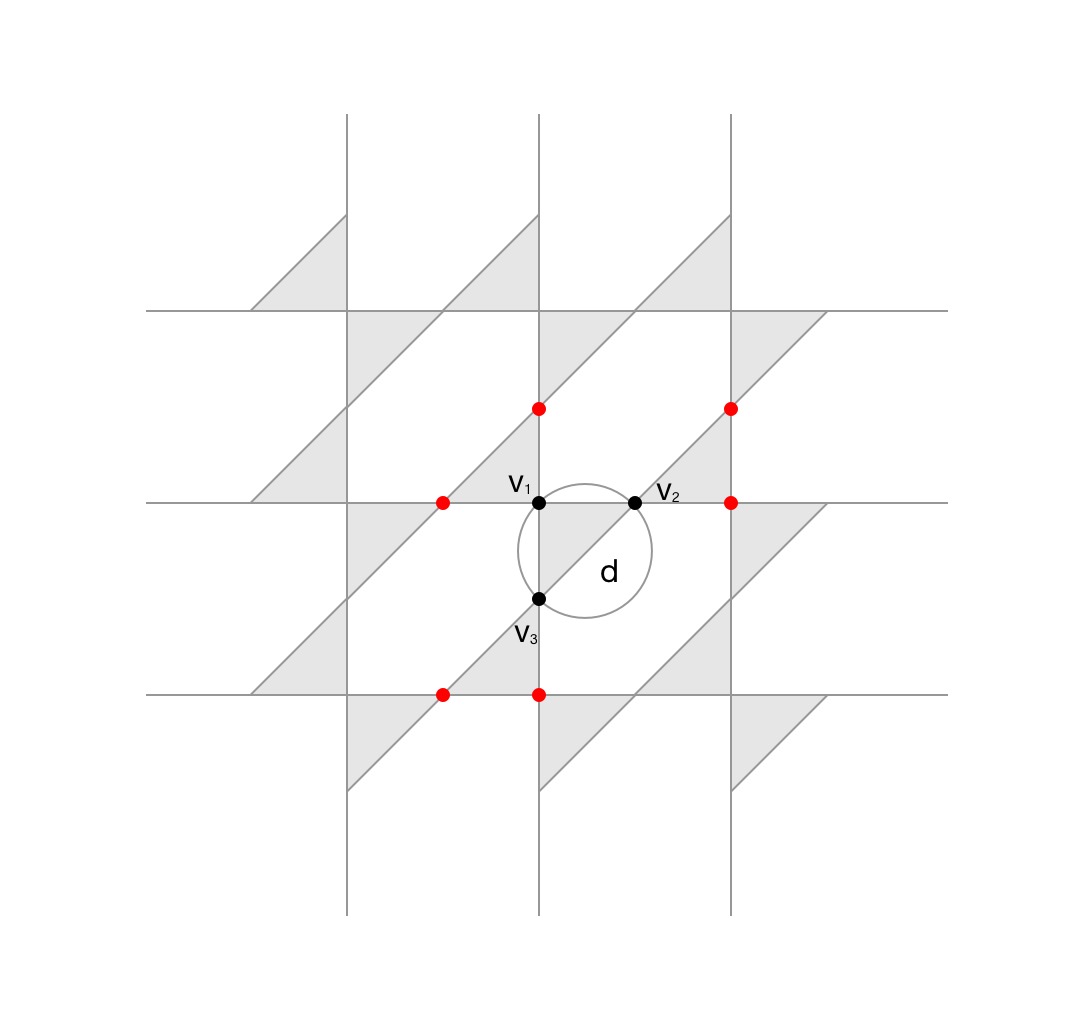}
 \caption{Example of $B(v_i,1)$ is the circle in black and the boundary of the union over all $i$ of $B(v_i,1)$ is colored in red.}
 \label{fig:boundball}
 \end{figure}

Let $d \in F_l^n$ and let $v_1,...,v_m$ be the vertices of $G_n$ which lie on 
the boundary of $d$, see Figure \ref{fig:boundball}.
Let $B(v,r) \subset G_n$ denote the ball centered at vertex $v$ of
radius $r$ in the path metric on $G_n$.  
It follows from the definition of $F_l^n$ and the fact that $G_n$ is the
planar dual of the graph of the disk pattern $G(D_n)$ (without the
vertex corresponding to the unbounded face), 
$d \in F_l^n$ implies that $B(v_i,l) \subset G_n$ but
$B(v_i, l+1) \not\subset G_n$ for $i=1,...,m$.  
Hence the distance from $v_i$ to
$\partial G_n$ is $l$ i.e. $v_i \in \partial B(x,l) $ for some
$x \in \partial G_n$ for all $i=1,...,m$. 


  Hence $F_l^n \subset \cup_{x \in \partial G_n} \partial B(x,l)$. 
%

\begin{lemma}
$$\displaystyle {\lim_{n \to \infty}} \frac{|\cup_l F_l^n|}{v(P_n)} = 1.$$ 
\end{lemma}

\begin{proof}
 We begin by showing that there exists $m>0$ such that $|\partial B(x,l)| \leq  ml$ for any $x \in G_n$. 
 By definition of F\o lner convergence,
  $G_n \subset \Gamma_{\mathcal{L}}$.  Babai \cite{Babai} showed that
    the growth rate for almost vertex transitive graphs with one end
    is quadratic i.e.  growth of $|B(x,l)|$ is quadratic in $l$. Since
    $\Gamma_{\mathcal{L}}$ is a biperiodic 4-valent planar graph, it
      satisfies the conditions of Babai's theorem, and hence has
      quadratic growth rate. By definition, the vertices in
      $\partial B(x,l)$ are incident to vertices in $B(x,l-1)$, hence
      $|\partial B(x,l)|$ has linear growth rate in $l$. 


Thus, $|F_l^n| \leq ml |\partial G_n |$ and
 we obtain

$$\displaystyle {\lim_{n \to \infty}} \frac{|F_l^n|}{v(P_n)} \leq \displaystyle {\lim_{n \to \infty}} \frac{ml |\partial G_n|}{3 a(K_n)} = \frac{ml}{3} \displaystyle {\lim_{n \to \infty}} \frac{|\partial G_n|}{|G_n|} \cdot \frac{|G_n|}{a (K_n)} = 0.$$  
 
Since $G_n \subset G(\mathcal{L})$, every vertex of $G_n$ lies on a disk in $F_l^n$ for some $l$ and for every disk in $F_l^n$ there are no vertices in $G(K_n) - G_n$ which lie on the disk. Now, by assumption $\displaystyle {\lim_{n \to \infty}} \frac{|G_n|}{3a(K_n)} = 1$. Hence $\displaystyle {\lim_{n \to \infty}} \frac{|\cup_l F_l^n|}{v(P_n)} = \displaystyle {\lim_{n \to \infty}} \frac{|G_n|}{3a(K_n)} = 1.$ 
\end{proof}


 

Let $f_k^n$ be the face with $k$ sides that is not contained in
$\cup_l F_l^n$ which does not meet the point at infinity. For each 
$n$, 
$\vol(C(f_k^n)) \leq k\lambda(\pi/6)$, where 
$\lambda(\theta)$ is the 
the Lobachevsky function defined as 
 $$\lambda(\theta) = - \int_{0}^{\theta} \log|2\sin(t)| dt,$$
 whose maximum value is $\lambda(\pi/6)$ \cite{Thurston} (also see \cite{CA1}).

Let $E^n$ denote the sum of the actual volumes of all the cones over the faces $f_k^n$, for every integer $k$. Then we have 
 
 \begin{equation}\label{eqn:zero}
E^n \leq \sum_k \sum_{f_k^n} k \lambda(\pi/6) = \sum_k k |f_k^n| \lambda(\pi/6).
 \end{equation}
 
 As mentioned before, every vertex of $G_n$ lies on a disk in $F_l^n$
 for some $l$ and for every disk in $F_l^n$ there are no vertices in
 $\Gamma_{K_n} - G_n$ which lie on the disk. By assumption
 $G_n \subset \Gamma_{\mathcal{L}} \cap (n\Lambda)$.
 where $n\Lambda$
 represents $n^2$ copies of the fundamental domain for the
 lattice $\Lambda$ such that $L = \mathcal{L}/\Lambda$.  

\begin{figure}
 \centering
 \includegraphics [height=5cm]{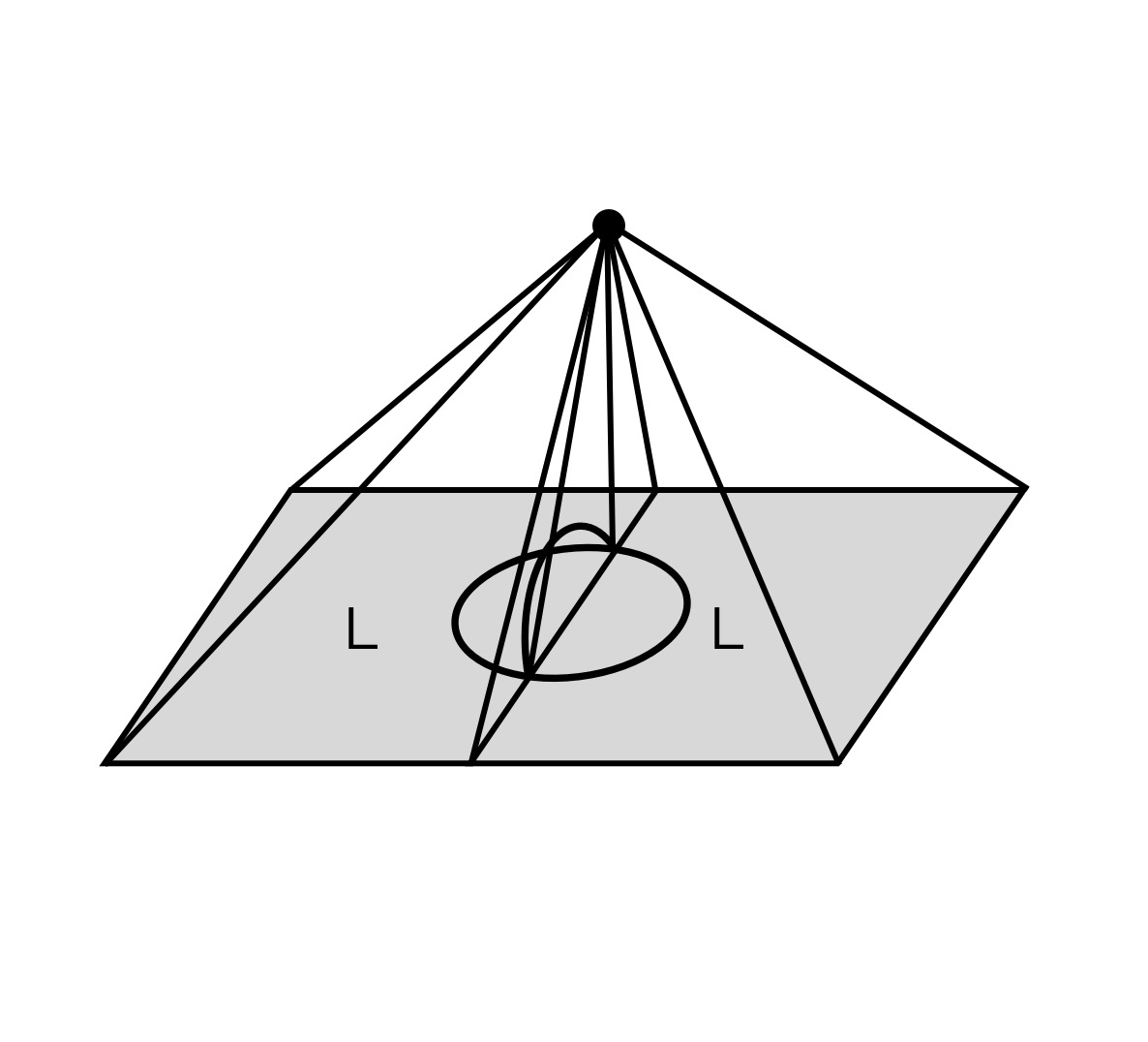}
 \caption{Two copies of the link $L$ coned to the point at infinity. The disk pattern from one copy of $L$ extends to the next copy.}
 \label{fig:2L}
 \end{figure}

 Since the cone vertex of the torihedron for $T^2 \times I - L$ is
 at infinity, the disk pattern obtained from taking $n^2$ copies
 of $L$
 just extends the disk pattern from one copy of $L$ to $n \times n$
 grid as in Figure \ref{fig:2L}. The graph for the disk pattern for
 $n^2$ copies of $L$ intersects $\Gamma_{K_n}$ in $G_n$, as in Figure
 \ref{fig:nbynL}.
 
 For any face $f$ in $F_l^n$, let $\delta_l^n$ be a positive number such that $\vol(C(f)) = \vol(C(f')) \pm \delta_l^n$ where $f'$ is a face in the disk pattern of $\mathcal{L}$ such that the graph isomorphism between $G(D_n)$ and $G(D_{\infty})$ sends $f$ to $f'$. Furthermore, we choose $\delta_l^n$ so that we can bound the sequence of $\delta_l^n$ by a sequence which will converge to zero as in Lemma \ref{lem:atkinson}.

Then  
\begin{equation}\label{eqn:two}
\vol (P_n) = \sum_l \sum_{f \in F_l^n} (\vol (f') \pm \delta_l^n) + E^n.
 \end{equation}
By Equation \ref{eqn:two} we get the following:

 \begin{equation}  \label{eqn:three}
  \vol (P_n) = \frac{1}{2} n^2 \vol((T^2 \times I)-L) + \sum_l \sum_{f \in F_l^n} (\pm \delta_l^n) + E^n.
 \end{equation}
 
 We divide each term by $a(K_n)$ and take the limit. For the first term of Equation \ref{eqn:three}, we obtain

 $$\displaystyle {\lim_{n \to \infty}} \frac{1}{2} \frac{n^2 \vol ((T^2 \times I)-L)}{a(K_n)} = \frac{1}{2} \frac{n^2 \vol((T^2 \times I)-L)}{n^2 a(L)} = \frac{1}{2} \frac{\vol((T^2 \times I)-L)}{a(L)} .$$

From our assumption of the F\o lner convergence the last condition gives us: 
$$ {\lim_{n \to \infty}}  \frac{a(K_n)}{n^2 a(L)} =1.  $$
 
 By Lemma \ref{lem:atkinson} there are positive numbers $\epsilon_l$ such that $\delta_l^n \leq \epsilon_l,$ so the second term of equation \ref{eqn:three} becomes 
 
 $$ \displaystyle {\lim_{n \to \infty}} \frac{|\sum_l \sum_{f \in F_l^n} (\pm \delta_l^n) |}{a(K_n)}  \leq \displaystyle {\lim_{n \to \infty}} \frac{\sum_l |F_l^n| \epsilon_l}{a(K_n)}.$$

\begin{lemma}

 $\displaystyle {\lim_{n \to \infty} \frac{\sum_l |F_l^n|
   \epsilon_l}{a(K_n)}} =0$. 
  \end{lemma}
 \begin{proof}
 Fix any $\epsilon >0.$ Because
 $\displaystyle {\lim_{l \to \infty}} \epsilon_l = 0,$ there is $K$
 sufficiently large such that $\epsilon_l < \epsilon/3$, for $l >
 K$.
 Then $\sum_{l=1}^K \epsilon_l$ is a finite number, say $M$. Since
 we've seen above that
 $\displaystyle {\lim_{n \to \infty}} \frac{\cup_l |F_l^n|}{v(P_n)} =
 1$
 and $\displaystyle {\lim_{n \to \infty}} \frac{|F_l^n|}{v(P_n)}=0$,
 there exist $N$ such that if $n>N$ then
 $ \displaystyle{{\rm max}_{l\leq L} \frac{|F_l^n|}{v(P_n)} < \frac{\epsilon}{(3M \cdot K)}}$ and
 $\displaystyle{\frac{|\cup_l F_l^n|}{v(P_n)}< (1+ \epsilon)}$. Then for $n>N$,
 
 $$\frac{\sum_l |F_l^n| \epsilon_l}{v(P_n)} = \frac{\sum_{l=1}^K |F_l^n| \epsilon_l}{v(P_n)} + \frac{\sum_{l>K}|F_l^n| \epsilon_l}{v(P_n)} < \frac{\epsilon K}{3M \cdot K} + (1+ \epsilon)\frac{\epsilon}{3}< \epsilon.$$
 
 \end{proof} 
 
Now setting $v(P_n) = 3a(K_n)$ we get that the limit of the second term is zero. 
\ \\
\indent Finally, by Equation \ref{eqn:one} and Equation \ref{eqn:zero} we get the third term of Equation \ref{eqn:three} equal to zero  $$\displaystyle {\lim_{n \to \infty}} \frac{E^n}{a(K_n)} \leq \displaystyle {\lim_{n \to \infty}} \frac{\sum_k k|f_k^n| \lambda(\pi/6)}{a(K_n)} =0.  $$ 

Therefore, $\displaystyle {\lim_{n \to \infty}} \frac{\vol(P_n)}{a(K_n)} = \frac{1}{2} \frac{\vol(T^2 \times I - L)}{a(L)}$ which means $\displaystyle {\lim_{n \to \infty}} \frac{\vol(K_n)}{a(K_n)} = \frac{\vol(T^2 \times I - L)}{a(L)}.$ 

\qed

Recall by $\mathcal{W}_f$ we mean the fully augmented square weave link whose quotient is $W_f$ with volume $10c\vtet$

\begin{corollary} \label{cor:falAsymptoticVolume}
Let $K_n$ be any sequence of hyperbolic fully augmented links such that $K_n$ F\o lner converges everywhere to $\mathcal{W}_f$. Then $$\displaystyle {\lim_{n \to \infty}}  \frac{\vol(K_n)}{a(K_n)} = 10\vtet.$$
\end{corollary}

{\it Proof.} The Corollary follows from Proposition \ref{prop:falSquareWeave} and Theorem \ref{thm:volDetConjFAL}.

\qed

\bibliographystyle{plain}
\bibliography{references-ak}

\end{document}